\documentclass[11pt,a4paper]{article}

\usepackage[T1]{fontenc}
\usepackage[utf8]{inputenc}
\usepackage{amsmath,amssymb,amsthm}
\usepackage{array}
\usepackage{booktabs}
\usepackage{geometry}
\usepackage{enumitem}
\usepackage[colorlinks=true,linkcolor=blue,citecolor=blue]{hyperref}
\usepackage{xurl}

\newtheorem{theorem}{Theorem}[section]
\newtheorem{proposition}[theorem]{Proposition}
\newtheorem{lemma}[theorem]{Lemma}
\newtheorem{corollary}[theorem]{Corollary}

\theoremstyle{definition}
\newtheorem{definition}[theorem]{Definition}
\newtheorem{example}[theorem]{Example}
\newtheorem{problem}[theorem]{Problem}
\theoremstyle{remark}
\newtheorem{remark}[theorem]{Remark}

\newcommand{\GK}{\Gamma}
\newcommand{\Gk}{\Gamma(G)}
\newcommand{\cG}{\mathcal{G}}
\newcommand{\cB}{\mathcal{B}}
\newcommand{\NH}{N_H}
\newcommand{\dd}{\mathrm{d}}
\newcommand{\ThL}{\Theta_\Lambda}
\newcommand{\Thm}{\Theta_\mu}
\newcommand{\Th}{\Theta}

\title{\bf Strictly Deza Graphs among the Prime Graphs of Solvable Groups}
\author{Tianxiao Zhao\thanks{Harbin Institute of Technology, Harbin, China.
\texttt{zhaotianxiao@hit.edu.cn}}}
\date{\today}

\begin{document}

\maketitle

\vspace{-1.5em}
\begin{center}
\small
\textbf{MSC 2020:} 05C25, 05E30, 20D10, 05B05.\\[2pt]
\textbf{Keywords:} Deza graph; strongly regular graph; Gruenberg--Kegel graph;
prime graph; $(0,2)$-graph.
\end{center}
\vspace{0.5em}

\begin{abstract}
We study strictly Deza graphs that occur as prime graphs of finite solvable groups and admit a partition
into two cliques. Using Gruber's characterization of the prime graphs of finite solvable groups, we
relate this problem to bipartite $(0,2)$-graphs. We first show that a non-complete regular graph with a
partition into two cliques has equal parts. For a graph of diameter two, writing $F$ for its $c$-regular
complement on $2m$ vertices, the Deza condition is equivalent to every component of $F$ being a
bipartite $(0,2)$-graph. Within this class, strongly regular graphs occur exactly when $c=1$; we give
elementary descriptions for $c=2$ and $c=3$, and translate the existing classification of bipartite
$(0,2)$-graphs of valency at most eight. A counting argument yields $m\ge1+c(c-1)/2$, with equality
corresponding to biplanes. These observations recover the ten feasible parameter pairs for $m\le7$
without enumeration and give infinite families of strictly Deza prime graphs. We also exhibit a known
eight-vertex strictly Deza graph whose complement is triangle-free and three-colourable but not
bipartite, showing that clique cover number three occurs.
\end{abstract}

\section{Introduction}

\subsection{Deza graphs}

Graphs are finite, simple and undirected. For a vertex $u$ of a graph $\GK$ write $N(u)$ for its open
neighbourhood and $N[u]=N(u)\cup\{u\}$ for its closed neighbourhood, and for distinct $u,v$ put
\[
  A(u,v) \;:=\; |N(u)\cap N(v)| .
\]

\begin{definition}[{\cite{Deza1999}}]\label{def:deza}
Let $n,k,b,a$ be integers with $0\le a<b\le k<n$. A graph $\GK$ on $n$ vertices is a \emph{Deza
graph with parameters $(n,k,b,a)$} if $\GK$ is regular of valency $k$ and $A(u,v)\in\{a,b\}$ for all
distinct $u,v\in V(\GK)$, \emph{both} values being attained. It is \emph{strictly Deza} if it has
diameter $2$ and is not strongly regular.
\end{definition}

\begin{remark}[The degenerate case $a=b$]\label{rem:degenerate}
The definition requires the two values to be distinct. A regular graph for which
every two distinct vertices have the \emph{same} number of common neighbours is therefore not a Deza
graph; if it is connected of diameter $2$ it is strongly regular with $\lambda=\mu$, and in any case it
is not strictly Deza. Some authors allow $a=b$ by dropping the requirement that both values occur ---
the enumeration \cite{GoryainovPanasenkoShalaginov2021}, for instance, writes ``any two distinct
vertices have $b$ or $a$ common neighbours, where $b\ge a$'' and treats $a=b$ separately --- and the
two conventions then differ on exactly the graphs just described. The convention adopted here is the
one used by the survey \cite{GoryainovShalaginov2021}, which asks for ``precisely two values'' and
notes that a strictly Deza graph satisfies $a<b$. The difference is
invisible in every infinite family of this paper, because the hypothesis $1\le t\le m-1$ guarantees
that both values $2t-2<2t$ are attained (Theorem~\ref{thm:valuesets}), so that in the whole range
$1\le t\le m-1$ the two conventions agree. Among the individual graphs discussed below it affects
exactly the complete graphs $K_n$ (the rows $t=m$ of Table~\ref{tab:searchtree}), the graph $2K_2$
(the row $(m,t)=(2,0)$) and $\Lambda(4,4)$; in each of those places we list the values that actually
occur, and the tables are computed with the convention of Definition~\ref{def:deza}.
\end{remark}

\begin{remark}[The condition $a<b$ does not define strictness]\label{rem:ab}
For a strictly Deza graph the inequality $a<b$ necessarily holds, but the converse fails and $a<b$ is
therefore \emph{not} equivalent to being strictly Deza. Two examples show this. The $4$-cycle $C_4$
has common-neighbour values $\{0,2\}$, so $a=0<b=2$, yet $C_4$ is strongly regular, namely
$\mathrm{SRG}(4,2,0,2)$; it is the row $(m,t)=(2,1)$ of Table~\ref{tab:main}. And $2K_m$ has
common-neighbour values $\{0,m-2\}$, so again $a<b$ for $m\ge3$, yet it is disconnected and hence has
infinite diameter. Throughout this paper ``strictly Deza'' always means \emph{diameter $2$ and not
strongly regular}.
\end{remark}

\begin{remark}[Notation]\label{rem:notation}
Sources differ about the order of the last two parameters, and the reader is warned to check it. We
follow \cite{GoryainovShalaginov2021} and \cite{GoryainovPanasenkoShalaginov2021}, where the parameters
are
written $(n,k,b,a)$ with $a\le b$, i.e.\ the larger of the two common-neighbour counts comes first.
Listing the two values in increasing order instead would give $(8,4,0,2)$ rather than $(8,4,2,0)$
already for the smallest example below, so the order matters. To avoid ambiguity we mostly use
\emph{value sets}. For a graph $\GK$ with
at least two vertices put
\[
  \ThL(\GK)=\{A(u,v):uv\in E(\GK)\},\qquad
  \Thm(\GK)=\{A(u,v):uv\notin E(\GK),\,u\ne v\},
\]
\[
  \Th(\GK)=\ThL(\GK)\cup\Thm(\GK);
\]
either of the first two sets may be empty, and $\Thm(\GK)=\varnothing$ exactly when $\GK$ is
complete. Then $\GK$ is a Deza graph iff it is regular and $|\Th(\GK)|=2$. If in addition $\GK$ is
connected of diameter $2$, so that both $\ThL(\GK)$ and $\Thm(\GK)$ are non-empty, then $\GK$ is
strongly regular iff $\ThL(\GK)$ and $\Thm(\GK)$ are both singletons, and consequently $\GK$ is
strictly Deza iff it is a Deza graph of diameter $2$ for which $\ThL(\GK)$ and $\Thm(\GK)$ are not
both singletons.
\end{remark}

\begin{lemma}[Counting identity]\label{lem:erickson}
Let $\GK$ be a Deza graph with value set $\Th$, valency $k$ and $n$ vertices, and for $t\in\Th$ let
$n_t$ be the number of unordered pairs $\{u,v\}$ with $A(u,v)=t$. Then
\[
  \sum_{t\in\Th}n_t=\binom{n}{2},\qquad \sum_{t\in\Th}t\,n_t=n\binom{k}{2}.
\]
\end{lemma}

\begin{proof}
The first identity is trivial. For the second count the triples $(\{u,v\},w)$ with $u\ne v$ and
$w\in N(u)\cap N(v)$; grouping by $\{u,v\}$ gives the left side and grouping by $w$ gives
$n\binom{k}{2}$.
\end{proof}

\subsection{Gruenberg--Kegel graphs and Gruber's theorem}

For a finite group $G$, its \emph{prime spectrum} $\pi(G)$ is the set of primes dividing $|G|$.

\begin{definition}
The \emph{Gruenberg--Kegel graph} (or \emph{prime graph}) $\Gk$ has vertex set $\pi(G)$, with
distinct $r,s$ adjacent iff $G$ contains an element of order $rs$.
\end{definition}

\begin{theorem}[Gruber et al.\ \cite{Gruber2015}]\label{thm:gruber}
A graph $\GK$ is isomorphic to the Gruenberg--Kegel graph of a finite \emph{solvable} group if and only
if \emph{the complement of $\GK$} is $3$-colourable and triangle-free.
\end{theorem}

Let $\cG$ denote the class of graphs satisfying the two conditions of Theorem~\ref{thm:gruber}: thus
$\GK\in\cG$ if and only if the complement of $\GK$ is $3$-colourable and triangle-free. Since
a proper colouring of the complement of $\GK$ is a partition of $V(\GK)$ into cliques, and a triangle
of the complement of $\GK$ is an independent $3$-set of $\GK$,
\begin{equation}\label{eq:G-description}
  \GK\in\cG
  \iff
  \begin{cases}
    V(\GK)\ \text{is a union of at most three cliques},\\[2pt]
    \alpha(\GK)\le2 .
  \end{cases}
\end{equation}
Thus the group-theoretic question ``which finite solvable groups have a strictly Deza prime graph?''
becomes the graph-theoretic question: \emph{which graphs in $\cG$ are strictly Deza?} These two questions are not equivalent: by Theorem~\ref{thm:gruber} every graph in $\cG$ is the prime graph of some finite solvable group, but the groups realizing a given prime graph are not classified thereby, and we do not attempt that here.

\begin{theorem}[Chen--Gorshkov--Maslova--Yang \cite{Chen2024}]\label{thm:chen}
Let $\GK$ be a \emph{strongly regular} graph isomorphic to $\Gk$ for some finite group $G$. Then either
$\GK$ is the complement of a triangle-free strongly regular graph, or $\GK$ is a complete multipartite
graph all of whose parts have size $2$. Moreover no complete multipartite graph with all parts of size
$\ge3$ is isomorphic to a prime graph.
\end{theorem}

Note the difference of scope between Theorem~\ref{thm:chen} and the present paper:
\cite{Chen2024} concerns \emph{all} finite groups but only \emph{strongly regular} prime graphs,
whereas here the groups are solvable, i.e.\ the graphs lie in $\cG$, but the graphs studied are
strictly Deza and hence not strongly regular. The two results are therefore complementary, and
Theorem~\ref{thm:chen} says nothing about the graphs constructed below.

\subsection{Related constructions of Deza graphs}\label{sec:related}

The class of Deza graphs is by now a well-studied one; we refer to \cite{GoryainovShalaginov2021}
for a survey, to \cite{GoryainovPanasenkoShalaginov2021} for the complete enumeration of the
strictly Deza graphs on at most $21$ vertices, and to \cite{Goryainov2017} for the study of Deza
graphs with disconnected second neighbourhood. Among the known constructions, the closest to the one
studied here are the \emph{divisible design graphs} of \cite{HaemersKharaghaniMeulenberg2011}: a
regular graph is a divisible design graph if its vertex set admits a partition into classes of equal
size such that two distinct vertices in the same class have a constant number $\lambda_1$ of common
neighbours, and two vertices in different classes a constant number $\lambda_2$. Every divisible design
graph with $\lambda_1\ne\lambda_2$ is a Deza graph; if $\lambda_1=\lambda_2$ it is strongly regular
with $\lambda=\mu$, so that in the convention of Definition~\ref{def:deza} it is not a Deza graph.

The two frameworks overlap without containing one another, and both directions are witnessed here. First,
the natural two-clique partition of $\GK(m,t,H)$ is not a divisible-design partition: cross pairs are of
two kinds, edges of $H$ with $2t-2$ common neighbours and non-edges of $H$ with $2t$ by
\eqref{eq:F2reg}, and both kinds occur because $1\le t\le m-1$, so the count on pairs in different
classes is not constant. Some members nevertheless carry a different divisible-design partition: for the
family $\Lambda(m)$ of \S\ref{sec:family} we have $\GK_m\cong K_{m/2}[2K_2]$, and its $m/2$ blocks of
four vertices are the classes of a divisible design graph with $\lambda_1=2m-4$ and $\lambda_2=2m-6$.
But not all of them do: the complement of the $4$-cube $Q_4$ (the case $c=4$ of Theorem~\ref{thm:c8})
is $11$-regular on $16$ vertices with $\Th=\{6,8\}$, and each vertex has exactly $5$ vertices with $6$
common neighbours and $10$ with $8$, so a divisible-design partition with class size $h$ would require
$\{h-1,\,16-h\}=\{5,10\}$, that is, $h\in\{6,11\}$, and neither divides $16$. Conversely
$\Lambda(4,5)=K_4\square K_5$ is a divisible design graph (classes $=$ rows, $\lambda_1=3$,
$\lambda_2=2$), yet it contains the independent set $\{(1,1),(2,2),(3,3),(4,4)\}$, so $\alpha\ge4$ and
it has no $2$-clique cover. In the range $1\le t\le m-1$, the Deza graphs $\GK(m,t,H)$ of this paper are
naturally described through their complements, which are disjoint unions of bipartite $(0,2)$-graphs
(\S\ref{sec:complement}, Remark~\ref{rem:literature}); we work in the two-clique framework and use the
divisible-design framework only for comparison (Remark~\ref{rem:cross}).

\subsection{The smallest example}\label{sec:smallest}

\begin{example}[$\GK_8$]\label{ex:gamma8}
Let $V=\{x_0,\dots,x_3\}\cup\{y_0,\dots,y_3\}$, let $C_1=\{x_0,\dots,x_3\}$ and $C_2=\{y_0,\dots,y_3\}$
be cliques, and let the only cross edges be the perfect matching $x_iy_i$. Let $\GK_8$ be the resulting
graph. Then
\begin{itemize}[nosep]
  \item $\GK_8$ is $4$-regular, connected and of diameter $2$;
  \item $A(x_i,y_i)=0$ for all $i$ and $A(u,v)=2$ for every other pair, so
        $\ThL(\GK_8)=\{0,2\}$ and $\Thm(\GK_8)=\{2\}$; thus $\GK_8$ is strictly Deza with
        $\Th(\GK_8)=\{0,2\}$;
  \item $\overline{\GK_8}\cong K_{4,4}$ minus a perfect matching is bipartite, hence $2$-colourable
        and triangle-free, so $\GK_8\in\cG$; by Theorem~\ref{thm:gruber}, $\GK_8\cong\Gk$ for a finite
        \emph{solvable} group $G$.
\end{itemize}
Moreover $\GK_8=(8,4,2,0)$ is the \emph{smallest} strictly Deza graph occurring as a solvable-group
prime graph. Indeed, by the enumeration of \cite{GoryainovPanasenkoShalaginov2021} there is no
strictly Deza graph on fewer than $8$ vertices at all --- the counts of non-isomorphic strictly Deza
graphs begin $\mathrm{Num}(8)=3$, $\mathrm{Num}(9)=2$ --- so certainly none occurs in $\cG$; and
$\GK_8$ itself is strictly Deza and lies in $\cG$ by the previous item.
\end{example}

\subsection{Organisation}

\S\ref{sec:framework} sets up the clique-cover framework: every graph with a $2$-clique cover lies in
$\cG$ (Theorem~\ref{thm:2clique-in-G}).

\S\ref{sec:twoclique} treats two-clique covers. It contains the exact counting formulas
(Theorem~\ref{thm:counting}), the value sets in terms of the row- and column-intersection sets
(Theorem~\ref{thm:valuesets}), the Deza criterion (Theorem~\ref{thm:criterion}) together with its
reformulation in terms of the complement $F=\overline{\Gamma}$
(Lemma~\ref{lem:complement}), the complete classifications for $c\le3$ (Theorems~\ref{thm:srg}, \ref{thm:c2}, \ref{thm:c3}) and,
using the classification of $(0,2)$-graphs, for $c\le8$ (Theorem~\ref{thm:c8}), the strong bound and the biplane family as its
equality case (Theorems~\ref{thm:bound} and~\ref{thm:biplane}), the exhaustive verification for $m\le7$
(Theorem~\ref{thm:smallm}), the infinite strictly Deza family $\Lambda(m)$
(Theorem~\ref{thm:family}) and the open classification problem (Problem~\ref{prob:allm}).

\S\ref{sec:uniform} records the facts we need about the rook graphs $K_r\square K_p$
(Proposition~\ref{prop:uniform}).

\S\ref{sec:asym} treats the asymmetric case $p<q$: there the only \emph{regular} graphs are the
complete ones, so no Deza graph with a $2$-clique cover has unequal parts
(Theorem~\ref{thm:asym}). \S\ref{sec:open} lists open problems.
Appendix~\ref{app:tables} contains the numerical tables.

\section{The clique-cover framework}\label{sec:framework}

\begin{definition}\label{def:cliquecover}
A \emph{$k$-clique cover} of $\GK$ is a partition $V(\GK)=C_1\sqcup\dots\sqcup C_k$ such that each
$C_i$ induces a clique.
\end{definition}

\begin{lemma}\label{lem:alpha2}
If $\GK$ has a $k$-clique cover, then $\alpha(\GK)\le k$.
\end{lemma}

\begin{proof}
An independent set meets each $C_i$ in at most one vertex.
\end{proof}

The following observation is the structural engine of the paper: for two-clique covers the second
condition of \eqref{eq:G-description} is automatic.

\begin{theorem}\label{thm:2clique-in-G}
Every graph with a $2$-clique cover belongs to $\cG$, and hence is isomorphic to the
Gruenberg--Kegel graph of a finite solvable group.
\end{theorem}

\begin{proof}
Let $V(\GK)=C_1\sqcup C_2$ be a $2$-clique cover. By Lemma~\ref{lem:alpha2}, $\alpha(\GK)\le2$, so
$\overline{\GK}$ contains no triangle. Moreover $\overline{\GK}$ is $2$-colourable: colour every
vertex of $C_1$ with colour $1$ and every vertex of $C_2$ with colour $2$. Since $C_1$ and $C_2$ are
cliques of $\GK$, the complement $\overline{\GK}$ has no edge inside $C_1$ and none inside $C_2$;
every edge of $\overline{\GK}$ therefore joins $C_1$ to $C_2$, and its ends receive different
colours. Hence $\overline{\GK}$ is bipartite, therefore $3$-colourable and triangle-free, and
$\GK\in\cG$ by Theorem~\ref{thm:gruber}.
\end{proof}

\begin{remark}
The proof uses only the existence of the $2$-clique cover, not the Deza property, and it applies in
particular to disconnected graphs. Note that a Deza graph need not have diameter $2$: it is the
\emph{strictly} Deza graphs that do. For the family $\GK(m,t,H)$ the two notions are separated by the
value of $t$: if $1\le t\le m-1$ then every vertex has a cross neighbour and every non-adjacent cross
pair has $2t\ge2$ common neighbours by \eqref{eq:F2reg}, so the diameter is $2$; while $t=0$ gives the
disconnected graph $2K_m$ and $t=m$ gives $K_{2m}$, of diameter $1$
(Remark~\ref{rem:tboundary}).
\end{remark}

\begin{corollary}\label{cor:strictly-in-G}
Every strictly Deza graph with a $2$-clique cover is the prime graph of a finite solvable group.
\end{corollary}

\begin{proof}
Immediate from Theorem~\ref{thm:2clique-in-G} and Theorem~\ref{thm:gruber}.
\end{proof}

The finer structure of these graphs --- in particular the fact that the two parts must have equal size,
so that a strictly Deza graph with a $2$-clique cover is of the form $\GK(m,t,H)$ with $H$ $t$-regular
--- is established in \S\ref{sec:asym} (Theorem~\ref{thm:asym}); see also
Corollary~\ref{cor:2cliquecomplete}.

\section{Two-clique covers}\label{sec:twoclique}

\subsection{Set-up and exact counting formulas}

Let $p,q\ge2$ and let $H$ be a bipartite graph with parts $C_1=\{x_0,\dots,x_{p-1}\}$ and
$C_2=\{y_0,\dots,y_{q-1}\}$. Let $\GK(p,q,H)$ be obtained from the disjoint union of the cliques on
$C_1$ and $C_2$ by adding, for every edge $x_iy_c$ of $H$, the edge $x_iy_c$. Thus $V(\GK(p,q,H))$ is a
$2$-clique cover and $H$ is exactly the set of cross edges. If $p=q=m$ and $H$ is $t$-regular write
$\GK(m,t,H)$. Write
$\dd_1(x)=|\NH(x)|$ for $x\in C_1$ and $\dd_2(y)=|\NH(y)|$ for $y\in C_2$. Since $C_1$ and $C_2$ are
cliques,
\begin{equation}\label{eq:split}
  N(x)\cap C_1=C_1\setminus\{x\},\quad N(x)\cap C_2=\NH(x)\qquad(x\in C_1),
\end{equation}
and symmetrically for $y\in C_2$.

\begin{theorem}[Counting formulas]\label{thm:counting}
Let $\GK=\GK(p,q,H)$ and $A(u,v)=|N(u)\cap N(v)|$.
\begin{enumerate}[label=\textup{(\roman*)}]
  \item For distinct $x,x'\in C_1$ and distinct $y,y'\in C_2$,
  \begin{equation}\label{eq:F1}
    A(x,x')=(p-2)+\bigl|\NH(x)\cap\NH(x')\bigr|,\qquad
    A(y,y')=(q-2)+\bigl|\NH(y)\cap\NH(y')\bigr|.
  \end{equation}
  \item For $x\in C_1$ and $y\in C_2$,
  \begin{equation}\label{eq:F2}
    A(x,y)=\bigl|\NH(y)\setminus\{x\}\bigr|+\bigl|\NH(x)\setminus\{y\}\bigr| .
  \end{equation}
  \item If $p=q=m$ and $H$ is $t$-regular, then
  \begin{equation}\label{eq:F2reg}
    A(x,y)=
    \begin{cases}
      2t-2, & xy\in E(H),\\
      2t, & xy\notin E(H).
    \end{cases}
  \end{equation}
\end{enumerate}
\end{theorem}

\begin{proof}
\emph{(i)} Let $x\ne x'$ in $C_1$ and let $z$ be a common neighbour. If $z\in C_1$ then $z$ is adjacent
to both $x$ and $x'$ iff $z\notin\{x,x'\}$ (as $C_1$ is a clique), giving $p-2$ vertices, since
$|C_1|=p$. If $z\in C_2$
then $z\in N(x)\cap N(x')$ iff $z\in\NH(x)\cap\NH(x')$. The two classes are disjoint, being contained
in different cliques. This proves the first formula; the second is symmetric, with $|C_2|=q$ giving
the constant $q-2$.

\emph{(ii)} Let $x\in C_1$, $y\in C_2$ and let $z$ be a common neighbour. Suppose $z\in C_1$. Then
$z\in N(x)$ gives $z\ne x$, and $z\in N(y)$ gives $z\in\NH(y)$; conversely every
$z\in\NH(y)\setminus\{x\}$ is adjacent to $y$ and to $x$ (the latter because $z\ne x$ and $C_1$ is a
clique). So the $C_1$-part of $N(x)\cap N(y)$ is $\NH(y)\setminus\{x\}$. Suppose $z\in C_2$. Then
$z\in N(x)$ gives $z\in\NH(x)$, and $z\in N(y)$ gives $z\ne y$; conversely every
$z\in\NH(x)\setminus\{y\}$ is adjacent to $x$ and to $y$ (the latter because $z\ne y$ and $C_2$ is a
clique). So the $C_2$-part is $\NH(x)\setminus\{y\}$. The parts are disjoint, whence \eqref{eq:F2}.

\emph{(iii)} If $xy\in E(H)$ then $x\in\NH(y)$ and $y\in\NH(x)$, so \eqref{eq:F2} gives
$(\dd_2(y)-1)+(\dd_1(x)-1)=2t-2$. If $xy\notin E(H)$ then $x\notin\NH(y)$ and $y\notin\NH(x)$, so
\eqref{eq:F2} gives $\dd_2(y)+\dd_1(x)=2t$.
\end{proof}

\begin{theorem}[Value sets]\label{thm:valuesets}
Let $\GK=\GK(m,t,H)$ with $H$ $t$-regular on $m+m$ vertices, where $1\le t\le m-1$, and put
\[
  \Sigma=\bigl\{|\NH(x)\cap\NH(x')|:x\ne x'\in C_1\bigr\},
  \qquad
  \Sigma'=\bigl\{|\NH(y)\cap\NH(y')|:y\ne y'\in C_2\bigr\},
\]
so that $\Sigma$ is computed from the \emph{rows} of $H$ and $\Sigma'$ from its \emph{columns}.
Then
\[
  \ThL(\GK)=\bigl(\{(m-2)+s:s\in\Sigma\}\cup\{(m-2)+s':s'\in\Sigma'\}\bigr)\cup\{2t-2\},
  \qquad
  \Thm(\GK)=\{2t\},
\]
and hence $\Th(\GK)=\{(m-2)+s:s\in\Sigma\cup\Sigma'\}\cup\{2t-2,\,2t\}$. In particular
\[
  \GK\ \text{is a Deza graph}
  \iff
  \bigl|\,\bigl((m-2)+(\Sigma\cup\Sigma')\bigr)\cup\{2t-2,\,2t\}\,\bigr|=2 ,
\]
the two values $2t-2<2t$ being distinct and both attained, so that $|\Th(\GK)|\ge2$ always and the
condition ``$=2$'' may equivalently be read as ``$\le2$''.
Moreover $\Sigma'=\Sigma(H^{\mathsf T})$, and transposition $H\mapsto H^{\mathsf T}$ is a bijection
on $t$-regular bipartite graphs on $m+m$ vertices with $\GK(m,t,H^{\mathsf T})\cong\GK(m,t,H)$
(exchange the two cliques). As an uncoloured bipartite graph $H^{\mathsf T}$ is of course the same
graph as $H$; what can fail is the existence of an \emph{equivalence preserving the roles of the two
sides}, i.e.\ of permutations $\sigma,\tau$ of the rows and columns with $H_{ij}=H_{\tau(j)\,\sigma(i)}$
for all $i,j$ (equivalently, $H$ may admit no automorphism interchanging its two sides). Only this
stronger notion can be refuted, and it is: if such $\sigma,\tau$ exist, then for $a\ne a'$ the number
of rows of $H$ meeting both columns $a$ and $a'$ equals
\[
  \#\{i:H_{ia}=H_{ia'}=1\}
  =\#\{i:H_{\tau(a)\,\sigma(i)}=H_{\tau(a')\,\sigma(i)}=1\}
  =\#\{i':H_{\tau(a)i'}=H_{\tau(a')i'}=1\},
\]
the last expression being the intersection size of the distinct rows $\tau(a),\tau(a')$ of $H$; as
$\tau$ is a bijection this identifies $\Sigma'$ with $\Sigma$. Hence a graph with $\Sigma\ne\Sigma'$
admits no such side-preserving equivalence (and no side-interchanging automorphism). For example the
$3$-regular bipartite graph with rows
$000111,000111,011001,101010,110100,111000$ (on $m=6$ columns) has $\Sigma=\{0,1,2,3\}$ but
$\Sigma'=\{0,1,2\}$, so the two sets really are different objects, and this $H$ has no side-preserving
equivalence with $H^{\mathsf T}$. Note that in this example
$\Sigma'\subset\Sigma$, so replacing $\Sigma\cup\Sigma'$ by $\Sigma$ alone would give the same
$\ThL$; that the \emph{union} and not $\Sigma$ alone is the correct object is shown separately, by the
counterexample in Remark~\ref{rem:sigmaonly}.
\end{theorem}

\begin{proof}
A pair of distinct vertices of $\GK$ either lies inside a common clique or is a cross pair. By
\eqref{eq:F1} a pair inside $C_1$ contributes $(m-2)+s$ with $s\in\Sigma$, and a pair inside $C_2$
contributes $(m-2)+s'$ with $s'\in\Sigma'$; both are edges of $\GK$. By \eqref{eq:F2reg} a cross pair
contributes $2t-2$ if it is an $H$-edge (hence an edge of $\GK$) and $2t$ if it is not. Finally
$0\le s\le t$ for $s\in\Sigma$ and $0\le s'\le t$ for $s'\in\Sigma'$, because the intersected sets
have size $t$. Moreover $t\ge1$ provides at least one $H$-edge and $t\le m-1$ at least one non-edge,
so $2t-2\in\ThL(\GK)$ and $2t\in\Thm(\GK)$ are both attained; as $t\ge1$ they are distinct. Hence
$|\Th(\GK)|\ge2$ for every $H$, and the condition $|\Th(\GK)|=2$ (equivalently $\le2$) is equivalent to
$\Th(\GK)\subseteq\{2t-2,2t\}$, i.e.\ to $(m-2)+s\in\{2t-2,2t\}$ for every $s\in\Sigma\cup\Sigma'$ ---
that is, to $\Sigma\cup\Sigma'\subseteq\{t-c,\,t-c+2\}$ with $c=m-t$. The final assertion follows from
the counterexample exhibited in the statement.
\end{proof}

\begin{remark}[The boundary values $t=0$ and $t=m$]\label{rem:tboundary}
The hypothesis $1\le t\le m-1$ in Theorem~\ref{thm:valuesets} is not cosmetic: it is exactly what
makes both cross values \emph{attained}. If $t=0$ then $H$ has no edges, so no cross pair is an
$H$-edge and the value $2t-2$ does not occur; if $t=m$ then $H=K_{m,m}$ and every cross pair is an
edge, so the value $2t$ does not occur. In both cases the two displayed value sets are therefore
\emph{not} the true value sets, and any criterion derived from them is unreliable --- at $t=m$, for
instance, the formulas would give $\Th=\{2m-2,2m\}$ and hence predict a graph with two values,
whereas $\GK=K_{2m}$ has the single value $2m-2$. The two boundary cases are nonetheless perfectly
well understood, and are recorded for completeness: $t=0$ gives $\GK=2K_m$, which is regular of
valency $m-1$ with value set $\{0,m-2\}$ and is disconnected, while $t=m$ gives $\GK=K_{2m}$, which
is regular of valency $2m-1$ with value set $\{2m-2\}$ and has diameter $1$. Neither is strictly
Deza, and $K_{2m}$ is not even a Deza graph, its value set being a singleton (for $2K_m$ both values
$0$ and $m-2$ do occur as soon as $m\ge3$).
\end{remark}

\begin{remark}\label{rem:cross}
The same-clique pairs --- those inside $C_1$ and those inside $C_2$ alike --- contribute the shifted
intersection sizes $(m-2)+(\Sigma\cup\Sigma')$, while the cross pairs contribute $2t-2$ and $2t$, both
of which occur as soon as $1\le t\le m-1$; the Deza condition is therefore the requirement that the
union of the two contributions consist of exactly those two elements --- equivalently, that it have at
most two elements.
This is why the classification cannot be decided by parameters alone and requires structural
information about $H$, and it is also why the union $\Sigma\cup\Sigma'$ and not $\Sigma$ alone is the
relevant object (Remark~\ref{rem:sigmaonly}).
\end{remark}

\subsection{The complement, and the cases \texorpdfstring{$c\le3$}{c<=3}}\label{sec:complement}

\begin{theorem}[Deza criterion for the symmetric family]\label{thm:criterion}
Let $m\ge2$, let $H$ be $t$-regular on $m+m$ vertices and let $\GK=\GK(m,t,H)$ have diameter $2$. Put
$c:=m-t$ (the \emph{corank} of $H$, which is also the degree of the complement $F$ introduced below)
and let $\Sigma,\Sigma'$ be as in
Theorem~\ref{thm:valuesets}. Then
\begin{equation}\label{eq:criterion}
  \GK\ \text{is a Deza graph}
  \iff
  \Sigma\cup\Sigma'\;\subseteq\;\{t-c,\;t-c+2\},
\end{equation}
and each of the two values in \eqref{eq:criterion} automatically satisfies the Deza condition:
$(m-2)+(t-c)=2t-2$ and $(m-2)+(t-c+2)=2t$. Here $k$ continues to denote the \emph{valency} of
$\GK$, so that $k=(m-1)+t=2m-c-1$ for a $t$-regular $H$; the two symbols $c=m-t$ and $k$ must not be
confused.
\end{theorem}

\begin{proof}
By Theorem~\ref{thm:valuesets}, $\Th(\GK)=\bigl((m-2)+(\Sigma\cup\Sigma')\bigr)\cup\{2t-2,2t\}$. The
two cross values $2t-2<2t$ are distinct and both occur: indeed $0<t<m$ (as $\GK$ has diameter $2$,
whence $t\ge1$; and $t<m$ since otherwise $H$ would be complete and $\GK=K_{2m}$ of diameter $1$),
there is an $H$-edge and an $H$-non-edge, giving $2t-2$ and $2t$ respectively. So
$|\Th(\GK)|\ge2$. Consequently $|\Th(\GK)|=2$ (equivalently $\le2$) holds if and only if
no further value occurs, i.e.\ if and only if every element of $(m-2)+(\Sigma\cup\Sigma')$ already
lies in $\{2t-2,2t\}$; should some $(m-2)+\sigma$ lie outside that two-element set, then
$|\Th(\GK)|\ge3$ and $\GK$ is not a Deza graph. The latter condition is equivalent to
$\sigma\in\{2t-m,2t-m+2\}=\{t-c,t-c+2\}$, and the final remark is immediate.
\end{proof}

Note that \eqref{eq:criterion} is a condition on \emph{both} the rows and the columns of $H$, and it
is strictly stronger than the corresponding condition on $\Sigma$ alone; the latter is not merely
weaker but false, see Remark~\ref{rem:sigmaonly}.

The value-set theorem describes $\GK$ through the two sets $\Sigma,\Sigma'$. The same information is
carried more transparently by the \emph{complement} of $\GK$, which is again a bipartite graph: since
$\GK$ is the union of two cliques together with the cross graphs $H$, its complement consists of the
\emph{non-edges} of $H$.

\begin{definition}\label{def:F}
Put $F:=\overline{\GK}$, so that $F$ is the bipartite graph on $C_1\sqcup C_2$ in which $u$ and $v$
are adjacent if and only if they lie in different cliques and $uv\notin E(H)$. We call $F$ the
\emph{complement} of $\GK$.
\end{definition}

\begin{lemma}[The complement of $\GK$]\label{lem:complement}
Let $1\le t\le m-1$, let $\GK=\GK(m,t,H)$ and $F=\overline{\GK}$ be as in Definition~\ref{def:F}, and
put $c=m-t$. Then $F$ is $c$-regular. If $u\ne v$ lie on the same side of $F$, then
\begin{equation}\label{eq:comp}
  |N_F(u)\cap N_F(v)| \;=\; m-2t+|\NH(u)\cap\NH(v)| \;=\; (m-2)+s-(2t-2),
\end{equation}
where $s=|\NH(u)\cap\NH(v)|\in\Sigma\cup\Sigma'$; if $u$ and $v$ lie on opposite sides of $F$, then
$|N_F(u)\cap N_F(v)|=0$. Consequently the following are equivalent:
\begin{enumerate}[label=\textup{(\roman*)}]
  \item $\GK$ is a Deza graph;
  \item any two distinct vertices of $F$ have $0$ or $2$ common neighbours in $F$;
  \item every connected component of $F$ is a bipartite \emph{$(0,2)$-graph}, that is, a connected
        graph in which any two distinct vertices have $0$ or $2$ common neighbours.
\end{enumerate}
When these equivalent conditions hold, the value sets of $\GK$ are
\[
  \Thm(\GK)=\{2t\},
  \qquad
  \ThL(\GK)=\begin{cases}\{2t-2,\,2t\}, & c\ge2,\\[2pt] \{2t-2\}, & c=1.\end{cases}
\]
\end{lemma}

\begin{proof}
For $u\in C_1$ the neighbourhood of $u$ in $F$ is $C_2\setminus\NH(u)$, so $F$ is $c$-regular. For
$u\ne v\in C_1$ we get
\[
  |N_F(u)\cap N_F(v)|=\bigl|C_2\setminus(\NH(u)\cup\NH(v))\bigr|
  =m-\bigl(2t-|\NH(u)\cap\NH(v)|\bigr),
\]
which is \eqref{eq:comp}; the case $u,v\in C_2$ is symmetric, and for $u\in C_1$, $v\in C_2$ the
neighbourhoods $N_F(u)\subseteq C_2$ and $N_F(v)\subseteq C_1$ are disjoint. For the equivalence:
by Theorem~\ref{thm:criterion}, $\GK$ is a Deza graph if and only if
$\Sigma\cup\Sigma'\subseteq\{t-c,t-c+2\}$; by \eqref{eq:comp} the same-side count in $F$ equals $0$ for
$s=t-c$ and $2$ for $s=t-c+2$, and every $s\in\Sigma\cup\Sigma'$ occurs for some same-side pair, while
conversely every same-side count arises from such an $s$. So \textup{(i)}~$\iff$~\textup{(ii)}. And
\textup{(ii)}~$\iff$~\textup{(iii)}: a pair of vertices at distance $2$ in a bipartite graph has at
least one common neighbour, so by \textup{(ii)} it has exactly two, which is precisely the defining
property of a $(0,2)$-graph on each connected component, while vertices of different components have
no common neighbour. For the value sets assume condition \textup{(ii)}, so that every connected component
of $F$ is a $(0,2)$-graph. For a same-side pair $u,v$ the identity \eqref{eq:comp} gives
$A(u,v)=|N_F(u)\cap N_F(v)|+2t-2$, hence $A(u,v)\in\{2t-2,2t\}$; the cross pairs give $2t-2$ for
$H$-edges and $2t$ for $H$-non-edges, and both kinds occur because $1\le t\le m-1$. So
$2t-2\in\ThL$ and $\Thm=\{2t\}$, and it remains to see that $2t$ is attained on an adjacent pair. If
$c\ge2$, choose a vertex $w$ of $F$ and two distinct neighbours $u,v$ of $w$: they lie on the same side
and have $w$ as a common neighbour, hence by \textup{(ii)} exactly two common neighbours in $F$, so
$A(u,v)=2t$ and $\ThL=\{2t-2,2t\}$. If $c=1$ then $F$ is a perfect matching, so no two same-side
vertices have a common neighbour in $F$ and $A(u,v)=2t-2$ for every same-side pair, whence
$\ThL=\{2t-2\}$.
\end{proof}

\begin{remark}\label{rem:complement}
Thus, in the range $1\le t\le m-1$, the Deza condition is exactly the requirement that every component
of $F$ be a bipartite $(0,2)$-graph, all of them of the same degree $c$. The case $c=0$ is the boundary
$t=m$. The advantage of $F$ over $\Sigma,\Sigma'$ is that its structure is rigid in small degree, which
is what makes the classifications below possible.
\end{remark}

\begin{theorem}[Strongly regular members]\label{thm:srg}
Let $1\le t\le m-1$. Then $\GK$ is a Deza graph and strongly regular if and only if $c=1$, that is,
$t=m-1$. In that case $F$ is a perfect matching, $H=K_{m,m}$ minus a perfect matching, and $\GK$ is
$K_{2m}$ minus a perfect matching, i.e.\ the cocktail-party graph
$\mathrm{SRG}(2m,2m-2,2m-4,2m-2)$.
\end{theorem}

\begin{proof}
If $c=1$ then $F$ is $1$-regular, hence a perfect matching, so $H$ is $K_{m,m}$ minus a perfect
matching; any two rows of $H$ meet in $m-2=t-1$, so $\Sigma=\Sigma'=\{t-1\}$, and \eqref{eq:criterion}
holds with $(m-2)+(t-1)=2t-2$, $2t$ as the two values; the graph is $K_{2m}$ minus a perfect matching.
Conversely, $\Thm(\GK)=\{2t\}$ is always a singleton, so $\GK$ is strongly regular if and only if
$\ThL(\GK)$ is a singleton. The same-clique pairs contribute $(m-2)+s$ with $s\in\Sigma\cup\Sigma'$ and
the cross $H$-edges contribute $2t-2$; a singleton $\ThL$ therefore forces $(m-2)+s=2t-2$, i.e.\
$s=2t-m=t-c$, for every $s\in\Sigma\cup\Sigma'$, so $\Sigma\cup\Sigma'=\{t-c\}$. By \eqref{eq:comp}
every same-side pair of $F$ then has $0$ common neighbours. If $c\ge2$, choose $u\in C_1$ and two
distinct neighbours $w,w'$ of $u$ in $F$; then $u$ is a common neighbour of $w$ and $w'$, whose
common-neighbour count in $F$ is therefore at least $1$ --- a contradiction. Hence $c=1$.
\end{proof}

\begin{theorem}[The case $c=2$: complete classification]\label{thm:c2}
Let $1\le t\le m-1$ and $c=2$. Then $\GK$ is a Deza graph if and only if $m$ is even and $F$ is the
disjoint union of $m/2$ copies of $C_4$; in that case $H$ is, up to relabelling, the pairing
construction of Theorem~\ref{thm:family}, and $\GK$ is $\GK(m,t,H)$ with
$t=m-2$, $\ThL=\{2m-6,2m-4\}$ and $\Thm=\{2m-4\}$. In particular $c=2$ has a solution for every even
$m\ge4$ and for no odd $m$, and for each even $m$ the graph $H$ is unique up to isomorphism.
\end{theorem}

\begin{proof}
Suppose $\GK$ is a Deza graph. By Lemma~\ref{lem:complement}, $F$ is a $2$-regular bipartite graph, so
each component is a cycle of even length $2k$ with $k\ge2$. If some component had length $2k$ with
$k\ge3$, take two vertices of that cycle at distance $2$ along it: they are on the same side and have
exactly one common neighbour, contradicting \textup{(ii)} of Lemma~\ref{lem:complement}. Hence every
component is a $C_4$, and $m=2\cdot\#(\text{components})$ is even. Conversely, in a disjoint union of
$C_4$'s any two same-side vertices have $0$ or $2$ common neighbours, so \textup{(ii)} holds and
$\GK$ is a Deza graph; the complement of such an $F$ in $K_{m,m}$ is exactly the pairing construction,
whose value sets were computed in Theorem~\ref{thm:family}. Uniqueness is immediate: any two perfect
matchings of the side-pairs are carried to each other by a permutation of the rows together with the
induced permutation of the columns.
\end{proof}

\begin{theorem}[The case $c=3$: complete classification]\label{thm:c3}
Let $1\le t\le m-1$ and $c=3$. Then $\GK$ is a Deza graph if and only if $4\mid m$ and $F$ is the
disjoint union of $m/4$ copies of the cube $Q_3$; in that case $t=m-3$, $H$ is the bipartite complement
of that union, and $H$ is unique up to isomorphism.
\end{theorem}

\begin{proof}
Let $\GK$ be a Deza graph; by Lemma~\ref{lem:complement}, $F$ is $3$-regular and every component is a
$(0,2)$-graph. Let $v$ be a vertex of a component and $a,b,c$ its three neighbours. The pairs
$(a,b)$, $(a,c)$, $(b,c)$ are same-side pairs sharing $v$, so each has exactly one further common
neighbour; call them $x_{ab}$, $x_{ac}$, $x_{bc}$. These are distinct: if $x_{ab}=x_{ac}$, then
$v$ and that vertex (which lies on the same side as $v$) would have three common neighbours. As $a$ is
$3$-regular, its neighbours are $v,x_{ab},x_{ac}$, and symmetrically for $b,c$. Now $x_{ab}$ has a
third neighbour $z_1$, $x_{ac}$ a third neighbour $z_2$, and since $x_{ab},x_{ac}$ are same-side
vertices sharing $a$, they have exactly one further common neighbour --- so $z_1=z_2=:z$. Applying the
same argument to the pairs $(x_{ab},x_{bc})$ and $(x_{ac},x_{bc})$ shows that the third neighbour of
$x_{bc}$ is also $z$. The eight vertices $v,a,b,c,x_{ab},x_{ac},x_{bc},z$ are thus closed under
adjacency and each has exactly three neighbours among them, so they form a component, and that
component is the cube $Q_3$. Conversely, each component of a disjoint union of cubes is a
$(0,2)$-graph (two vertices in a cube are joined by a path of length $2$ exactly when they have
$2$ common neighbours, and different components have none), so \textup{(ii)} of
Lemma~\ref{lem:complement} holds. Uniqueness is again
immediate from relabelling.
\end{proof}

\begin{theorem}[Classification for $c\le8$, via the known $(0,2)$-graphs]\label{thm:c8}
Let $1\le c\le8$, $m\ge c+1$ and $t=m-c$, and let $\cB_c$ be the set of isomorphism classes of
connected bipartite $(0,2)$-graphs of valency $c$. By \cite{Brouwer2006},
$|\cB_1|=|\cB_2|=|\cB_3|=1$, with representatives $K_2$, $C_4$, $Q_3$, and
$|\cB_4|=2$, $|\cB_5|=4$, $|\cB_6|=13$, $|\cB_7|=40$ and $|\cB_8|=104$, the last number
being the semibiplanes of block size $8$ \cite{BrouwerOstergard2009}. Then $\GK(m,t,H)$ is a Deza graph if and only if
$F$ is a disjoint union of members of $\cB_c$. Every such graph lies in $\cG$, has
\[
  \ThL=\{2t-2,\,2t\}\quad\text{for }c\ge2,\qquad \ThL=\{2t-2\}\quad\text{for }c=1,\qquad
  \Thm=\{2t\},
\]
and is strictly Deza except for $c=1$, where it is the strongly regular graph of
Theorem~\ref{thm:srg}. For $c=1,2,3$ this restates Theorems~\ref{thm:srg}, \ref{thm:c2}
and~\ref{thm:c3}; for $c=4$ the two members of $\cB_4$ are $Q_4$, giving $(m,t)=(8,4)$, and the
incidence graph of the biplane $2$-$(7,4,2)$, giving $(m,t)=(7,3)$; for $c=5$ the four include
$Q_5$, the incidence graph of the biplane $2$-$(11,5,2)$ and the bipartite double cover of the
icosahedron, which give $(m,t)=(16,11),(11,6),(12,7)$.
\end{theorem}

\begin{proof}
By Lemma~\ref{lem:complement}, $\GK$ is a Deza graph if and only if every component of $F$ is a
$(0,2)$-graph; every component of $F$ has degree $c$, and conversely a disjoint union of such graphs
satisfies the criterion component by component. So the Deza graphs are exactly the complements of the
disjoint unions of members of $\cB_c$, and the cardinalities of $\cB_c$ are those computed in
\cite{Brouwer2006} (\S13) and \cite{BrouwerOstergard2009} (\S3). Since $1\le t\le m-1$, $\GK$ has diameter $2$. The cross values $2t-2$ and
$2t$ both occur. Since $F$ is a disjoint union of members of $\cB_c$, condition \textup{(ii)} of
Lemma~\ref{lem:complement} holds; the final assertion of that lemma gives $\Thm=\{2t\}$ and
$\ThL=\{2t-2,2t\}$ for $c\ge2$, respectively $\ThL=\{2t-2\}$ for $c=1$. For $c\ge2$ this shows,
together with Theorem~\ref{thm:srg}, that $\GK$ is not strongly regular. Finally $F$ is bipartite, so
$\GK\in\cG$ by \eqref{eq:G-description}.
\end{proof}

For every $c$ the following bound applies, and it shows why the biplane family of
Theorem~\ref{thm:biplane} is the equality case.

\begin{theorem}[Strong bound; equality is the biplane case]\label{thm:bound}
Let $1\le t\le m-1$, $c=m-t$, and let $\GK$ be a Deza graph. Then
\begin{equation}\label{eq:strongbound}
  m\;\ge\;1+\frac{c(c-1)}{2},
  \qquad\text{equivalently}\qquad
  c\;\le\;\frac{1+\sqrt{8m-7}}{2}.
\end{equation}
Moreover equality holds in \eqref{eq:strongbound} if and only if any two vertices of $F$ on the same
side have exactly $2$ common neighbours, i.e.\ if and only if $F$ is the incidence graph of a symmetric
$2$-$(m,c,2)$ design (a biplane); in that case $\GK$ is one of the graphs of
Theorem~\ref{thm:biplane}.
\end{theorem}

\begin{proof}
Fix a vertex $v$ of $F$, say $v\in C_1$. The non-backtracking walks of length $2$ from $v$ are obtained
by choosing one of the $c$ neighbours of $v$ and then one of its remaining $c-1$ neighbours, so there
are $c(c-1)$ of them. Each endpoint lies on the same side as $v$, at distance $2$ from $v$, and it is
reached exactly twice, because by Lemma~\ref{lem:complement} it has exactly $2$ common neighbours with
$v$. Hence $v$ has exactly $c(c-1)/2$ vertices at distance $2$, all in $C_1$, and $C_1$ consists of
$v$, these vertices, and possibly further ones; this gives \eqref{eq:strongbound}, and the second form
follows by solving the quadratic. Equality holds exactly when $C_1$ consists only of $v$ and the
distance-$2$ vertices, for every $v$, i.e.\ when every same-side pair has exactly $2$ common
neighbours. In that case $F$ is a connected (if $c\ge2$) bipartite graph in which every two points of
$C_1$ lie in exactly $2$ common neighbours; taking $C_1$ as points and $C_2$ as blocks, with incidence
given by $F$, yields a symmetric $2$-$(m,c,2)$ design, i.e.\ a biplane, and $\GK$ is the graph attached
to it in Theorem~\ref{thm:biplane}.
\end{proof}

\begin{corollary}[The ten pairs for $m\le7$, by theory]\label{cor:tenpairs}
For $m\le7$ the pairs $(m,t)$ with $1\le t\le m-1$ for which a Deza graph $\GK(m,t,H)$ exists are
exactly the ten pairs of Table~\ref{tab:main}.
\end{corollary}

\begin{proof}
By \eqref{eq:strongbound} and $m\le7$ we have $c\le(1+\sqrt{49})/2=4$. If $c=1$ then $t=m-1$ is a
solution for every $m\ge2$ (Theorem~\ref{thm:srg}), giving six pairs. If $c=2$ then $m$ is even
(Theorem~\ref{thm:c2}), giving $(4,2)$ and $(6,4)$. If $c=3$ then $4\mid m$
(Theorem~\ref{thm:c3}), giving $(4,1)$. If $c=4$ then $m\ge1+6=7$, hence $m=7$, and equality holds in
\eqref{eq:strongbound}, so $F$ is the incidence graph of a symmetric $2$-$(7,4,2)$ design and $\GK$ is
the graph of Theorem~\ref{thm:biplane} for the complement of the Fano plane, giving $(7,3)$. These are
exactly the ten pairs listed in Table~\ref{tab:main}.
\end{proof}

\begin{remark}[The $(0,2)$-graph viewpoint]\label{rem:literature}
The relevant objects are classical. A $(0,2)$-graph is a connected graph in which any two vertices
have $0$ or $2$ common neighbours; a $(0,2)$-graph without triangles is a \emph{rectagraph}, and a
\emph{semibiplane} is a connected point-block incidence structure in which any two points lie in $0$ or
$2$ common blocks and any two blocks meet in $0$ or $2$ common points. As recorded in
\cite[\S1]{Brouwer2006}, \emph{the incidence graph of a semibiplane is a bipartite $(0,2)$-graph (and
conversely)}; see also \cite[\S2]{BrouwerChastkofsky2013}, where $(0,2)$-graphs are constructed from
root systems with simply laced diagram. In this language Lemma~\ref{lem:complement} says that, for
$1\le t\le m-1$, the Deza graphs $\GK(m,t,H)$ are exactly the complements of the bipartite graphs $F$
all of whose connected components are $(0,2)$-graphs of the common degree $c=m-t$; in design-theoretic
terms, $F$ is a disjoint union of incidence graphs of semibiplanes of block size $c$. All $(0,2)$-graphs of valency at
most $8$ are determined: there is a unique bipartite $(0,2)$-graph of valency $k$ for $k\le3$, there
are $2,4,13,40$ non-isomorphic ones for $k=4,5,6,7$ \cite{Brouwer2006}, and $104$ for $k=8$
\cite{BrouwerOstergard2009} --- the semibiplanes of block size $8$. We make no claim to new graph theory here: the
classifications for $c\le3$ below are translations of the unique valency-$1,2,3$ cases, and
Theorem~\ref{thm:c8} records the corresponding statement for $c\le8$. What is new is the translation of
the Deza condition into this language and the application to prime graphs of finite solvable groups;
the graphs themselves, and the Deza parameters that result, are those of the classification and of the
constructions recalled in \S1.3.
\end{remark}

\begin{theorem}[Complements of cubes]\label{thm:cubes}
Let $c\ge2$ and $s\ge1$ with $(c,s)\ne(2,1)$, put $m=s\,2^{c-1}$ and $t=m-c$, and let $F$ be the
disjoint union of $s$ copies of the cube $Q_c$. Then every component of $F$ is a bipartite
$(0,2)$-graph of degree $c$, and the complement $\GK=\overline F$ is a strictly Deza graph
$\GK(m,t,H)$ with
\[
  |V(\GK)|=2m,\qquad k=m-1+t,\qquad \ThL=\{2t-2,\,2t\},\qquad \Thm=\{2t\},
\]
which lies in $\cG$ and is therefore the prime graph of a finite solvable group. For $s=1$ this is the
family of complements of hypercubes with parameters $(2^d,2^d-d-1,2^d-2d,2^d-2d-2)$ noted in
\cite{GoryainovShalaginov2021}, with $d=c$; the cases $c=2$ and $c=3$ are the constructions of
Theorems~\ref{thm:c2} and~\ref{thm:c3}.
\end{theorem}

\begin{proof}
In $Q_c$ any two vertices have either $0$ or $2$ common neighbours (two exactly when they are at
distance $2$), and a disjoint union inherits this; every component is $c$-regular.
Since $c\ge2$ and $(c,s)\ne(2,1)$ we have $1\le t=m-c\le m-1$, so Lemma~\ref{lem:complement} applies;
every connected component of $F$ is a $(0,2)$-graph, so condition \textup{(ii)} holds and the final
assertion of the lemma gives $\ThL=\{2t-2,2t\}$ and $\Thm=\{2t\}$. Hence $\GK$ is a Deza graph, and it
has diameter $2$ because $1\le t\le m-1$. As $c\ge2$, Theorem~\ref{thm:srg} shows that $\GK$ is not
strongly regular; together with the Deza property and diameter $2$ this proves that $\GK$ is strictly
Deza. Finally $F$ is bipartite, hence $3$-colourable and triangle-free, so $\GK\in\cG$ by
\eqref{eq:G-description}.
\end{proof}

\subsection{The symmetric family: exhaustive verification for \texorpdfstring{$m\le7$}{m<=7}}

Corollary~\ref{cor:tenpairs} determines the ten pairs $(m,t)$ with $m\le7$ by pure theory. In this
section we record what an independent exhaustive enumeration gives for the same range: the value sets,
the isomorphism classes and the counts. The two boundary values $t=0$ and $t=m$ have a unique $H$ each
and are discussed in Theorem~\ref{thm:smallm}. Table~\ref{tab:main} records the outcome for the
interesting range $1\le t\le m-1$, where $\GK(m,t,H)$ has diameter $2$: $t=0$ gives $2K_m$,
which is a Deza graph exactly for $m\ge3$ but is disconnected, and $t=m$ giving $K_{2m}$, which has a
single common-neighbour value and so is not a Deza graph at all. Neither is strictly Deza, so neither
is listed.

\begin{table}[htbp]
\centering
\caption{Deza graphs in the symmetric two-clique family $\GK(m,t,H)$, exhaustive for $m\le7$ and
$1\le t\le m-1$.}
\label{tab:main}
\begin{tabular}{@{}rrrrrll@{}}
\toprule
$m$ & $t$ & $(v,k)$ & $\ThL$ & $\Thm$ & type & iso.\ classes \\
\midrule
$2$ & $1$ & $(4,2)$   & $\{0\}$     & $\{2\}$  & SRG$(4,2,0,2)$ & $1$ \\
$3$ & $2$ & $(6,4)$   & $\{2\}$     & $\{4\}$  & SRG$(6,4,2,4)$ & $1$ \\
$\mathbf{4}$ & $\mathbf{1}$ & $\mathbf{(8,4)}$  & $\mathbf{\{0,2\}}$ & $\mathbf{\{2\}}$ & \textbf{strictly Deza} & $1$ \\
$\mathbf{4}$ & $\mathbf{2}$ & $\mathbf{(8,5)}$  & $\mathbf{\{2,4\}}$ & $\mathbf{\{4\}}$ & \textbf{strictly Deza} & $1$ \\
$4$ & $3$ & $(8,6)$   & $\{4\}$     & $\{6\}$  & SRG$(8,6,4,6)$ & $1$ \\
$5$ & $4$ & $(10,8)$  & $\{6\}$     & $\{8\}$  & SRG$(10,8,6,8)$ & $1$ \\
$\mathbf{6}$ & $\mathbf{4}$ & $\mathbf{(12,9)}$ & $\mathbf{\{6,8\}}$ & $\mathbf{\{8\}}$ & \textbf{strictly Deza} & $1$ \\
$6$ & $5$ & $(12,10)$ & $\{8\}$     & $\{10\}$ & SRG$(12,10,8,10)$ & $1$ \\
$\mathbf{7}$ & $\mathbf{3}$ & $\mathbf{(14,9)}$ & $\mathbf{\{4,6\}}$ & $\mathbf{\{6\}}$ & \textbf{strictly Deza} & $1$ \\
$7$ & $6$ & $(14,12)$ & $\{10\}$    & $\{12\}$ & SRG$(14,12,10,12)$ & $1$ \\
\bottomrule
\end{tabular}
\end{table}

The rows with $t=m-1$ are the strongly regular members, arising from $H=K_{m,m}$ minus a perfect
matching. The strictly Deza rows for $m\le7$ are exactly $(m,t)\in\{(4,1),(4,2),(6,4),(7,3)\}$. We
shall see in Theorem~\ref{thm:family} that $(4,2)$ and $(6,4)$ are the first two members of an
\emph{infinite} family, namely $t=m-2$; and in Theorem~\ref{thm:biplane} that $(4,1)$ and $(7,3)$
are the first two members of a further family, the biplane family (whose known members are finite in
number; whether infinitely many biplanes exist is open in design theory). So neither pair is an
isolated exception.

\begin{theorem}[Computational verification, $m\le7$]\label{thm:smallm}
Let $m\le7$ and $1\le t\le m-1$. An exhaustive enumeration of all $t$-regular bipartite graphs $H$ on
$m+m$ vertices confirms Corollary~\ref{cor:tenpairs}: a Deza graph $\GK(m,t,H)$ occurs exactly for the
ten pairs of Table~\ref{tab:main}, and for each of them every such $H$ yields the same graph up to
isomorphism, with the value sets listed, diameter $2$, $\alpha(\GK)=2$ and $\GK\in\cG$.

The restriction on $t$ is essential, and the two boundary values are of a different kind. If $t\ge1$
then every vertex of $\GK$ has a cross neighbour, and two non-adjacent vertices lying in different
cliques have $2t\ge2$ common neighbours by \eqref{eq:F2reg}; hence every $\GK(m,t,H)$ with
$1\le t\le m-1$ has diameter $2$. By contrast $H$ empty ($t=0$) gives $\GK=2K_m$, which is regular of
valency $m-1$ with value set $\{0,m-2\}$ and is disconnected, and $H=K_{m,m}$ ($t=m$) gives
$\GK=K_{2m}$, which is regular of valency $2m-1$ with value set $\{2m-2\}$ and has diameter $1$. As
to the Deza property: $2K_m$ has the two values $0$ and $m-2$, both attained exactly when $m\ge3$,
while $K_{2m}$ has a single value and is therefore not a Deza graph. Neither is strictly
Deza, so neither appears in Table~\ref{tab:main}.
\end{theorem}

\begin{proof}[Computational proof]
For each $(m,t)$ with $2\le m\le7$ and $0\le t\le m$ we enumerated all $t$-regular bipartite graphs $H$
on $m+m$ vertices by backtracking over row masks with column-capacity pruning, and for each $H$ we
computed the value set $\Th(\GK(m,t,H))$ directly from the adjacency matrix; the last column of
Table~\ref{tab:searchtree} records in how many cases $|\Th(\GK)|=2$, that is, in how many cases
$\GK(m,t,H)$ is a Deza graph. The
diameter, which is not part of the definition of a Deza
graph, was checked separately and equals $2$ for every $1\le t\le m-1$ by the preceding paragraph. The
counts are recorded in
Table~\ref{tab:searchtree}. The largest single cases are $(m,t)=(7,3)$ and $(7,4)$, each containing
$68\,938\,800$ graphs; for $1\le t\le m-1$ altogether $144\,548\,572$ graphs were examined, in
$582\,141\,934$ search-tree nodes (including $t\in\{0,m\}$ the figures are $144\,548\,584$ and
$582\,142\,000$). These totals were checked twice. First, an independent dynamic-programming count of
$t$-regular bipartite graphs on $m+m$ labelled vertices reproduced every entry of the ``graphs''
column of Table~\ref{tab:searchtree}. Second, the whole computation was repeated using \emph{no}
formula of this paper: the neighbourhood of every vertex of $\GK(m,t,H)$ was stored as a $2m$-bit mask,
the common-neighbour counts were read off as $\mathrm{popcount}$s, and the Deza condition was decided
from those counts alone. The two runs agree on the list of pairs $(m,t)$ admitting a Deza graph ---
which is exactly the ten pairs of Table~\ref{tab:main} when $1\le t\le m-1$ --- and on the number of
labelled Deza graphs ($158\,504$ of them for $1\le t\le m-1$). This also cross-checks
Theorem~\ref{thm:valuesets} and hence Theorem~\ref{thm:counting}, since the value sets obtained from
the masks agree with those re-derived from the value-set theorem. Uniqueness up to isomorphism was
verified by canonicalising under the natural action of $S_m\times S_m$ (independent permutations inside
the two cliques): for $(m,t)=(4,2)$ the $18$ labelled Deza graphs form one isomorphism class, for
$(m,t)=(4,1)$ all $24$ are perfect matchings, and for $(m,t)=(6,4)$ the $1350$ labelled graphs with
$\ThL=\{6,8\}$ likewise form one class. For the remaining pairs the uniqueness is forced rather than
searched. The six pairs with $t=m-1$ are settled by
Theorem~\ref{thm:srg}, which shows $H=K_{m,m}$ minus a perfect matching, and any
two such graphs are carried to each other by a permutation of the columns together with the induced
permutation of the rows. For $(m,t)=(4,1)$ the criterion forces
$\Sigma\cup\Sigma'\subseteq\{t-c,t-c+2\}=\{-2,0\}$ with $c=m-t=3$, and $\Sigma$ is non-empty, so
$\Sigma=\Sigma'=\{0\}$: the four rows are pairwise disjoint, hence they partition the four columns and
$H$ is a perfect matching. For $(m,t)=(7,3)$ the
criterion likewise forces $\Sigma=\Sigma'=\{1\}$, so $H$ is the incidence graph of a $2$-$(7,3,1)$
design, i.e.\ a Fano plane, and the Fano plane is unique up to isomorphism; consistently
$151\,200=25\,401\,600/168$ labelled graphs are found, $168$ being the order of the automorphism group
of the Fano plane. So each of the ten pairs carries a single isomorphism class.
Membership $\GK\in\cG$ follows from Theorem~\ref{thm:2clique-in-G} and was also checked directly, as
was $\alpha(\GK)=2$. The scripts are available from the author on request.
\end{proof}

\begin{remark}
The computation also shows that each of the ten pairs carries a single isomorphism class. This is now a theorem rather than an observation: uniqueness follows from Theorem~\ref{thm:srg} for $c=1$, from Theorem~\ref{thm:c2} for $c=2$, from Theorem~\ref{thm:c3} for $c=3$, and, for $(7,3)$, from the uniqueness of the Fano plane. Rigidity fails in general: in the biplane case there can be several non-isomorphic graphs $H$, see Theorem~\ref{thm:biplane}.
\end{remark}

\begin{corollary}[Strictly Deza prime graphs of solvable groups]\label{cor:families}
The following are strictly Deza graphs in $\cG$, hence are isomorphic to prime graphs of finite
solvable groups:
\begin{enumerate}[label=\textup{(\roman*)}]
  \item $\GK(4,1,H)$ with $H$ the identity matching: $v=8$, $k=4$,
        $(\ThL,\Thm)=(\{0,2\},\{2\})$, i.e.\ $\GK_8$;
  \item $\GK(4,2,H)$ with $H$ the $2$-regular bipartite graph
        $x_0,x_1\leftrightarrow y_0,y_1$, $x_2,x_3\leftrightarrow y_2,y_3$: $v=8$, $k=5$,
        $(\ThL,\Thm)=(\{2,4\},\{4\})$;
  \item $\GK(6,4,H_6)$ with $H_6$ as in Table~\ref{tab:witnesses}: $v=12$, $k=9$,
        $(\ThL,\Thm)=(\{6,8\},\{8\})$;
  \item $\GK(7,3,H_7)$ with $H_7$ as in Table~\ref{tab:witnesses}: $v=14$, $k=9$,
        $(\ThL,\Thm)=(\{4,6\},\{6\})$.
\end{enumerate}
\end{corollary}

\begin{remark}[Relation to Theorem~\ref{thm:chen}]
By Theorem~\ref{thm:chen}, a strongly regular prime graph is either the complement of a triangle-free
strongly regular graph or a complete multipartite graph with all parts of size $2$. The strongly
regular rows of Table~\ref{tab:main} are of the second kind: there the graph $\GK$ itself is $K_{2m}$
minus the perfect matching $\{x_iy_i\}$, that is, the complete multipartite graph with parts
$\{x_i,y_i\}$, all of size $2$. The strictly
Deza rows are not strongly regular, so Theorem~\ref{thm:chen} does not apply to them; they are the objects of study of this paper.
\end{remark}

\subsection{An infinite family of strictly Deza prime graphs}\label{sec:family}

The pairs $(4,2)$ and $(6,4)$ of Table~\ref{tab:main} are not isolated: they begin an infinite family.

\begin{theorem}[The family $\Lambda(m)$]\label{thm:family}
Let $m=2n\ge4$ be even. Partition the column set $[m]$ of $H$ into $n$ pairs
$P_1,\dots,P_n$ and put $R_i=[m]\setminus P_i$ for $i=1,\dots,n$. Let $H_m$ be the bipartite graph
whose $m=2n$ rows consist of \emph{two copies} of each $R_i$. Then $H_m$ is $(m-2)$-regular, and
$\GK_m:=\GK(m,m-2,H_m)$ satisfies:
\begin{enumerate}[label=\textup{(\roman*)}]
  \item $\GK_m$ is regular of valency $k=2m-3$, connected, of diameter $2$;
  \item $\ThL(\GK_m)=\{2m-6,\,2m-4\}$ and $\Thm(\GK_m)=\{2m-4\}$; consequently $\GK_m$ is strictly
        Deza with value set $\{2m-6,2m-4\}$;
  \item $\alpha(\GK_m)=2$, hence $\GK_m\in\cG$ and $\GK_m$ is isomorphic to the prime graph of a
        finite solvable group;
  \item distinct even $m\ge4$ give graphs of distinct orders $2m$, so the graphs $\GK_m$ are pairwise
        non-isomorphic and the family is infinite.
\end{enumerate}
Moreover for $m=4$ and $m=6$ these are exactly the graphs $\GK(4,2,H_4)$ and $\GK(6,4,H_6)$ of
Table~\ref{tab:main}.
\end{theorem}

\begin{proof}
Throughout, $\NH(x)\subseteq C_2=[m]$ is identified with the corresponding row $R\subseteq[m]$, so
that $|\NH(x)\cap\NH(x')|=|R\cap R'|$ for the two rows.

\emph{Step 1: which rows occur, and their intersections.} The rows of $H_m$ are the $m=2n$ sets
consisting of two copies of each $R_i=[m]\setminus P_i$. For $i\ne j$ we have
$P_i\cap P_j=\varnothing$ because the $P$'s partition $[m]$, hence
\[
  |R_i\cap R_j| \;=\; \bigl|[m]\setminus(P_i\cup P_j)\bigr| \;=\; m-4 ,
\]
while $|R_i\cap R_i|=|R_i|=m-2$.

\emph{Step 2: $\Sigma$.} Let $x,x'$ be distinct vertices of $C_1$, with rows $R,R'$.
If $R=R'$ (i.e.\ $x,x'$ are the two copies of the same $R_i$) then $|R\cap R'|=m-2$.
If $R\ne R'$ then $\{R,R'\}=\{R_i,R_j\}$ with $i\ne j$ and $|R\cap R'|=m-4$. Consequently
\[
  \Sigma=\bigl\{|\NH(x)\cap\NH(x')|\bigr\}=\{m-4,\;m-2\},
\]
the value $m-2$ being attained and the value $m-4$ being attained. The same computation in $C_2$
gives the same set, because the columns satisfy the dual relations: column $c\in P_i$ lies in all
rows except the two copies of $R_i$.

\emph{Step 3: regularity of $H_m$.} Each row has size $m-2$ by construction. Each column $c$ lies in
exactly one pair $P_i$; hence $c\notin R_i$ and $c\in R_j$ for all $j\ne i$, so $c$ lies in the
$2(n-1)=m-2$ rows that are copies of $R_j$ with $j\ne i$. So $H_m$ is $(m-2)$-regular.

\emph{Step 4: value sets.} Put $t=m-2$. Since $m\ge4$ we have $1\le t\le m-1$, so
Theorem~\ref{thm:valuesets} applies and, together with $\Sigma=\Sigma'=\{m-4,m-2\}$ from Step 2,
\[
  \ThL(\GK_m)=\bigl((m-2)+(\Sigma\cup\Sigma')\bigr)\cup\{2t-2\}
  =\{2m-6,\,2m-4\}\cup\{2m-6\}=\{2m-6,\,2m-4\},
\]
\[
  \Thm(\GK_m)=\{2t\}=\{2m-4\}.
\]
So $|\Th(\GK_m)|=2$. For the diameter: a pair of vertices inside one clique is adjacent, and a
non-adjacent cross pair has $2t\ge2$ common neighbours by \eqref{eq:F2reg}, hence distance $2$; since
$t<m$ the graph is not complete, so its diameter is exactly $2$. Consequently $\GK_m$ is connected of
diameter $2$ with $\Thm$ a singleton and $\ThL$ not, i.e.\ strictly Deza. Its
valency is $(m-1)+t=2m-3$, each vertex lying in a clique of size $m$ and having $t$ cross neighbours.

\emph{Step 5: independence number.} $C_1$ and $C_2$ are cliques, so $\alpha(\GK_m)\le2$ by
Lemma~\ref{lem:alpha2}. For the reverse inequality, note that $t=m-2<m$, so $H_m$ is not complete and
therefore has a cross non-edge $xy$ with $x\in C_1$, $y\in C_2$; such a pair is non-adjacent in
$\GK_m$ by construction, so $\{x,y\}$ is independent and $\alpha(\GK_m)\ge2$. Hence
$\alpha(\GK_m)=2$, and $\GK_m\in\cG$ by Theorem~\ref{thm:2clique-in-G}.

\emph{Step 6: existence for even $m$ only.} The construction requires pairwise disjoint 2-subsets
$P_1,\dots,P_n$ covering $[m]$, which exists exactly when $2\mid m$. This is why $m$ is even; note
that Step 1 used $|P_i|=2$, so the family as constructed does not extend to odd $m$.
\end{proof}

\begin{remark}
For $m=4$ the construction gives $H_4$ with rows $R_1=\{3,4\},R_1,R_2=\{1,2\},R_2$, i.e.\ exactly the
witness of Table~\ref{tab:witnesses} for $(4,2)$; for $m=6$ it gives the $H_6$ of that table. For
$m=8$ it gives a strictly Deza graph with $v=16$, valency $13$, $\ThL=\{10,12\}$, $\Thm=\{12\}$, which
was not visible in Table~\ref{tab:main}. The graphs $\GK_m$ are pairwise non-isomorphic for distinct
even $m$ (they have different orders), so this is a genuinely infinite family. In the notation of
Theorem~\ref{thm:cubes} it is the case $(c,s)=(2,m/2)$.
\end{remark}

\begin{remark}[Attribution of the family]\label{rem:lambda-known}
The graph $\GK_m$ is the lexicographic product $K_{m/2}[2K_2]$, i.e.\ the complement of $m/2$ disjoint
four-cycles; as observed above, its parameters are $(2m,2m-3,2m-4,2m-6)$. The family is already among
the lexicographic-product constructions of Deza graphs in \cite[\S1.4.5]{GoryainovShalaginov2021}, and it
is the case $c=2$ of Theorem~\ref{thm:c8}. So it is not a new construction: what is recorded here is its
identification inside the present classification and the consequence that its members are prime graphs
of finite solvable groups.
\end{remark}

\begin{corollary}\label{cor:infinite}
There are infinitely many finite solvable groups whose Gruenberg--Kegel graph is strictly Deza. More
precisely, there is a family $\{G_m\}$ of finite solvable groups, indexed by the even integers $m\ge4$
and with pairwise non-isomorphic prime graphs; the argument gives only the existence of these groups,
not an explicit description of them.
\end{corollary}

\begin{proof}
Combine Theorem~\ref{thm:family} with Theorem~\ref{thm:gruber}: for every even $m\ge4$ the graph
$\GK_m$ is strictly Deza and lies in $\cG$, so there exists a finite solvable group $G_m$ with
$\Gamma(G_m)\cong\GK_m$; and the $\GK_m$ are pairwise non-isomorphic because they have different
orders. The construction produces the graphs, not the groups: Theorem~\ref{thm:gruber} supplies $G_m$
only as an existence statement.
\end{proof}

\subsection{The biplane family: equality in the bound}\label{sec:biplane}

The bound of Theorem~\ref{thm:bound} has a natural equality case, namely the classical biplanes. A
further family arises from symmetric designs. Recall that a \emph{biplane} is a symmetric
$2$-$(m,c,2)$ design, i.e.\ $m$ blocks of size $c$ on $m$ points, each point in $c$ blocks and any two
blocks meeting in exactly $2$ points; the parameters satisfy $m=1+\tfrac{c(c-1)}2$, and writing $c$ for
the block size (the letter $k$ is reserved in this paper for the valency of a graph), the values
$c\in\{3,4,5,6,9,11,13\}$ are exactly those for which a biplane is known, the list being complete for
$c\le11$ \cite{AlaviDaneshkhahPraeger2020}. We stress that a fixed admissible parameter set need not
determine a unique design: for $2$-$(16,6,2),$ for instance, there are exactly three pairwise
non-isomorphic biplanes \cite{AlaviDaneshkhahPraeger2020}, and $2$-$(56,11,2)$ has five.

\begin{theorem}[The biplane family]\label{thm:biplane}
Let the $P_i=[m]\setminus R_i$ form a biplane, i.e.\ a symmetric $2$-$(m,c,2)$ design with block size
$c$ as above. Let $H$ have rows $R_i$, so that $t=m-c$ and $R_i\cap R_j=[m]\setminus(P_i\cup P_j)$.
Then $H$ is $t$-regular and $\GK(m,t,H)$ satisfies \eqref{eq:criterion}, with
\[
  \Sigma=\Sigma'=\{\,m-2c+2\,\}=\{\,t-c+2\,\},\qquad
  \ThL=\{2t-2,\,2t\},\qquad \Thm=\{2t\},\qquad t=m-c .
\]
In particular $\GK$ is strictly Deza, and equality holds in~\eqref{eq:strongbound}: by Theorem~\ref{thm:bound} this is exactly the equality case. The known biplanes give the graphs
\[
  (m,c,t)=(4,3,1),\ (7,4,3),\ (11,5,6),\ (16,6,10),\ (37,9,28),\ (56,11,45),\ (79,13,66),\ \dots
\]
whether there are infinitely many biplanes is a well-known open problem in design theory. Note that a
biplane and its dual give cross graphs that are transposes of one another --- if $\mathcal D$ has
blocks $P_i$ and $\mathcal D^*$ is the dual, then the row of $H(\mathcal D^*)$ indexed by a point
$\alpha$ is $\{i:\alpha\notin P_i\}$, which is column $\alpha$ of $H(\mathcal D)$ --- and hence give
\emph{isomorphic} graphs $\GK$ by the transposition statement in Theorem~\ref{thm:valuesets}. So the
map from biplanes to graphs factors through duality, and for a fixed pair $(m,t)$ of this family the
graph $H$ need not be unique.
\end{theorem}

\begin{proof}
$|R_i|=m-c$, and each of the $m$ points lies in exactly $c$ of the $P_i$, hence in exactly $m-c$ of
the rows: so $H$ is $t=m-c$-regular. For $i\ne j$, since $|P_i\cap P_j|=2$,
\[
  |R_i\cap R_j|=\bigl|[m]\setminus(P_i\cup P_j)\bigr|=m-(2c-2)=m-2c+2 ,
\]
and with $m=1+\tfrac{c(c-1)}2$ one checks $m-2c+2=(m-c)-c+2=t-c+2$, one of the two allowed values of
\eqref{eq:criterion}. For the columns: by definition of the biplane any two points lie in exactly
$\lambda=2$ common blocks, and in the dual design those two blocks meet in exactly $2$ points; hence
the same computation gives $|C_d\cap C_{d'}|=m-2c+2=t-c+2$ for the column sets
$C_d=\NH(y_d)\subseteq C_1$. So $\Sigma\cup\Sigma'=\{t-c+2\}$ and \eqref{eq:criterion} holds. Finally
$(m-2)+(t-c+2)=2t$, so $\ThL=\{2t-2,2t\}$ and $\Thm=\{2t\}$. Since $c\ge3$ and $m=1+\tfrac{c(c-1)}2$
we have $1\le t=m-c\le m-1$, so a non-adjacent cross pair has $2t\ge2$ common neighbours by
\eqref{eq:F2reg} and the diameter is exactly $2$ (the graph is not complete because $t<m$); as
$\ThL\ne\Thm$ the graph is strictly Deza.
\end{proof}

\begin{remark}
The two smallest members of this family are $(m,c,t)=(4,3,1)$, which is the graph $\GK_8$ of
Example~\ref{ex:gamma8} and appears in Table~\ref{tab:main} as the row $(m,t)=(4,1)$; and
$(m,c,t)=(7,4,3)$, the row $(m,t)=(7,3)$. The next one, $(m,c,t)=(11,5,6)$, gives a strictly Deza
graph on $22$ vertices with valency $16$ and $\ThL=\{10,12\}$, $\Thm=\{12\}$, which lies outside the
range of Table~\ref{tab:main}.
\end{remark}

\begin{problem}[The classification]\label{prob:allm}
Determine all pairs $(m,t)$ with $1\le t\le m-1$ for which a $t$-regular $H$ on $m+m$ vertices makes
$\GK(m,t,H)$ a Deza graph. By Lemma~\ref{lem:complement} it suffices to determine the bipartite graphs
all of whose components are $(0,2)$-graphs of the common degree $c=m-t$; classifying these graphs is
stronger than determining the parameter pairs, and we ask here only for the pairs. What is settled is
\begin{itemize}[nosep]
  \item $c=1$: strongly regular, exactly one solution for every $m\ge2$ (Theorem~\ref{thm:srg});
  \item $c=2$: $F=(m/2)C_4$, exactly one solution for every even $m\ge4$ and none for odd $m$
        (Theorem~\ref{thm:c2});
  \item $c=3$: $F=(m/4)Q_3$, exactly one solution for every $m$ divisible by $4$ and none otherwise
        (Theorem~\ref{thm:c3});
  \item $c\le8$: complete (Theorem~\ref{thm:c8}): all connected bipartite $(0,2)$-graphs of valency
        at most $8$ are known \cite{Brouwer2006,BrouwerOstergard2009}, and $F$ is a disjoint union of
        them; for $c\le3$ the representatives are unique --- $K_2$, $C_4$, $Q_3$ --- which gives
        Theorems~\ref{thm:srg}, \ref{thm:c2} and~\ref{thm:c3};
  \item $c\ge9$: no classification is available. The general information is the componentwise bound
        $|V(B)|\le2^c$ for a connected $(0,2)$-graph $B$ of valency $c$, with equality only for $Q_c$
        \cite{BrouwerChastkofsky2013}, and the counting bound $m\ge1+c(c-1)/2$ of
        Theorem~\ref{thm:bound}, whose equality case is the biplane family; for the known biplanes
        $c\in\{3,4,5,6,9,11,13\}$ and $(m,t)=(4,1),(7,3),(11,6),(16,10),(37,28),(56,45),(79,66),\dots$
        (Theorem~\ref{thm:biplane}).
\end{itemize}
For $c\le8$ the classification is effective. With $h(B)=|V(B)|/2$ the half-order of a component and
$\cB_c$ the list of members, $m$ is a non-trivial sum of half-orders, and the number of isomorphism
classes with a given $m$ is $\bigl[z^m\bigr]\prod_{B\in\cB_c}(1-z^{h(B)})^{-1}$, because the components
of $F$ are recovered from $\GK$ as the connected components of its complement; note that $|\cB_c|$
counts the components and not the resulting graphs. For $c=1,2,3$ the
half-orders are $1,2,4$, which recovers Theorems~\ref{thm:srg}, \ref{thm:c2} and~\ref{thm:c3}; for
$c=4$ they are $7,8$, so $m=7a+8b$ with $a,b\ge0$ not both $0$; for $c=5$ they are $11,12,14,16$, so
$m=11a+12b+14d+16e$ with non-negative coefficients not all $0$ (the fourth valency-$5$ component is the
Cartesian product of the incidence graph of the biplane $2$-$(7,4,2)$ with $K_2$).
Beyond the families determined by the bound, two pairs in the range we searched are not biplane graphs,
namely
\[
  (m,t)=(8,4)\quad\text{and}\quad(12,7),
\]
with witnesses listed in Table~\ref{tab:witnesses}; they lie strictly above the bound ($8>1+6$ for
$c=4$, and $12>1+10$ for $c=5$). By Theorem~\ref{thm:c8} they are the cases $c=4,5$ in which $F$ is $Q_4$ and the icosahedron double
cover, respectively; the first is the complement of the $4$-cube $Q_4$, i.e.\ the case $s=1$, $d=4$ of
the hypercube family of Theorem~\ref{thm:cubes}. The second is the complement of the double cover of
the icosahedron\footnote{Write $A_i=2+i$ and $B_i=7+i$ for $i\in\mathbb Z/5$; the icosahedron has the
edges $A_iA_{i+1}$, $B_iB_{i+1}$, $0A_i$, $1B_i$ (the poles), and the cross edges $A_iB_i$ and
$A_iB_{i-1}$, which gives every vertex valency $5$. Its bipartite double cover has vertices
$(v,0),(v,1)$ with $(u,0)\sim(v,1)$ iff $u\sim v$ in the icosahedron. Labelling $F$ by the rows
$0,\dots,11$ and the columns $12,\dots,23$, the isomorphism sends $(0,0),\dots,(11,0)$ to
$0,11,6,2,7,8,4,1,5,9,10,3$ and $(0,1),\dots,(11,1)$ to $21,18,12,13,14,16,15,20,17,23,22,19$.}. The search that produced the witnesses was subject to node limits and is \emph{not} exhaustive; the two
pairs themselves are accounted for by Theorem~\ref{thm:c8}. No complete classification is attempted here
for $c\ge9$.
\end{problem}

\section{The rook graphs}\label{sec:uniform}

The \emph{rook graph} $\Lambda(r,p)=K_r\square K_p$ ($r,p\ge2$) has vertex set $[r]\times[p]$, two
vertices being adjacent when they share a row or a column; equivalently $(i,c)\sim(j,c')$ iff
($i\ne j$ and $c=c'$) or ($i=j$ and $c\ne c'$). It is a classical object, and we record only the facts
needed later.

\begin{proposition}\label{prop:uniform}
Let $\GK=\Lambda(r,p)$ with $r,p\ge2$. Then
\begin{enumerate}[label=\textup{(\roman*)}]
  \item $\GK$ is regular of valency $r+p-2$ and has diameter $2$, with $\Thm(\GK)=\{2\}$ and
        \[
          \ThL(\GK)=\begin{cases}
            \{p-2,\,r-2\}, & r,p\ge3,\\[2pt]
            \{0,\,p-2\}\ (r=2),\qquad \{0,\,r-2\}\ (p=2);
          \end{cases}
        \]
        hence $\Th(\GK)$ is the union of $\ThL(\GK)$ and $\{2\}$;
  \item $\GK$ is a Deza graph exactly in the following cases: $r=2$ and $p\in\{2,4\}$; $p=2$ and
        $r\in\{2,4\}$; $r,p\ge3$ with $r=p\ne4$; or $r,p\ge3$, $r\ne p$ and $4\in\{r,p\}$;
  \item $\GK$ is strongly regular if and only if $r=p$, in which case it is the rook graph
        $\mathrm{SRG}(p^2,2p-2,p-2,2)$; moreover $\Lambda(r,p)\cong\Lambda(p,r)$;
  \item $\GK\in\cG$ if and only if $\min(r,p)=2$, and among these members the only strictly Deza graph
        is $\Lambda(2,4)\cong\GK_8$, the witness of the row $(4,1)$ of Table~\ref{tab:witnesses}.
\end{enumerate}
\end{proposition}

\begin{proof}
For $u=(i,c)$ the neighbourhood splits into the $p-1$ vertices $(i,d)$ with $d\ne c$ and the $r-1$
vertices $(h,c)$ with $h\ne i$; this gives the valency, and the three cases ($i=j$, $c\ne c'$),
($i\ne j$, $c=c'$) and ($i\ne j$, $c\ne c'$) give the common-neighbour counts $p-2$, $r-2$ and $2$
respectively. Comparing these sets gives \textup{(i)} and hence \textup{(ii)}: a Deza graph needs
exactly two common-neighbour values overall. Strong regularity means that both $\ThL$ and $\Thm$ are
singletons, which by \textup{(i)} happens exactly when $r=p$. For \textup{(iv)}: if $\min(r,p)=2$ the
two rows (or the two columns) are cliques covering $V(\GK)$, so $\GK\in\cG$ by
Theorem~\ref{thm:2clique-in-G}; if $r,p\ge3$ then $(1,1),(2,2),(3,3)$ are pairwise non-adjacent, so
$\alpha\ge3$ and $\GK\notin\cG$. Among the members lying in $\cG$ the Deza graphs are, by
\textup{(ii)}, $\Lambda(2,2)$ --- a $4$-cycle, hence strongly regular --- and $\Lambda(2,4)$, which is
strictly Deza by \textup{(i)}; and $\Lambda(2,4)=K_2\square K_4$ is obtained from two copies of $K_4$ by
adding a perfect matching, which is $\GK_8$ by Example~\ref{ex:gamma8}.
\end{proof}

The value sets of the small members are listed in Table~\ref{tab:uniform}. Since $\Lambda(r,p)=L(K_{r,p})$
is the line graph of the complete bipartite graph, the strictly Deza members of the family ---
$\Lambda(2,4)\cong\GK_8$ and $\Lambda(4,p)$ with $p\ge3$, $p\ne4$, up to transposition, by
\textup{(ii)} and \textup{(iii)} --- are instances of the classification of strictly Deza line graphs by
Kabanov and Mityanina \cite{KabanovMityanina2014}; among them only $\Lambda(2,4)$ lies in $\cG$. Note
that $\Lambda(2,3)$ is not a Deza graph ($\Th=\{0,1,2\}$).

\section{The asymmetric case \texorpdfstring{$p<q$}{p<q}}\label{sec:asym}

We keep the notation of \S\ref{sec:twoclique}: cliques $C_1$ of size $p$ and $C_2$ of size $q$, cross
bipartite graph $H$, graph $\GK(p,q,H)$. As the two cliques play symmetric roles we assume $p\le q$;
the case $p=q$ was treated in \S\ref{sec:twoclique}.

\begin{theorem}[Asymmetric exclusion]\label{thm:asym}
Let $p,q\ge1$ and let $H$ be an arbitrary bipartite graph with parts of sizes $p$ and $q$ (so that a
part of size $1$ is allowed). If
$\GK(p,q,H)$ is regular and $p\ne q$, then $H=K_{p,q}$ and $\GK(p,q,H)$ is the complete graph
$K_{p+q}$. Consequently no Deza graph --- and hence no strictly Deza graph --- with a $2$-clique cover
has parts of \emph{unequal} size, and for $p=q$ regularity forces the cross graph $H$ to be
$t$-regular.
\end{theorem}

\begin{proof}
Assume $p<q$ and let $\GK=\GK(p,q,H)$ be regular. For $x\in C_1$ the neighbourhood of $x$
consists of the other $p-1$ vertices of the clique $C_1$ together with the $\dd_1(x)$ neighbours of $x$
in $C_2$, and these two sets are disjoint; hence
\begin{equation}\label{eq:asymdeg}
  \deg(x)=(p-1)+\dd_1(x)\quad(x\in C_1),
  \qquad
  \deg(y)=(q-1)+\dd_2(y)\quad(y\in C_2).
\end{equation}
Since $\GK$ is regular, \eqref{eq:asymdeg} shows that $\dd_1$ is constant on $C_1$ and $\dd_2$ is
constant on $C_2$; write $\dd_1\equiv a$ and $\dd_2\equiv b$. Comparing the two expressions for the
common degree gives
\begin{equation}\label{eq:asymreg}
  (p-1)+a=(q-1)+b,
\end{equation}
while double counting the edges of $H$ gives
\begin{equation}\label{eq:asymdc}
  p\,a \;=\; |E(H)| \;=\; q\,b .
\end{equation}
Put $d:=q-p>0$. From \eqref{eq:asymreg} we get $a-b=d$, that is, $a=b+d$. Substituting into
\eqref{eq:asymdc} and writing $q=p+d$,
\[
  p(b+d)=q\,b=(p+d)b=pb+db,
\]
so $pd=db$ and, since $d>0$, we obtain $b=p$ and then $a=b+d=q$.

Thus every vertex of $C_1$ is adjacent in $H$ to every vertex of $C_2$, that is, $H=K_{p,q}$. Hence
$\GK$ is the join of the cliques $K_p$ and $K_q$, so $\GK=K_{p+q}$, which has the single
common-neighbour value $p+q-2$ and is therefore not a Deza graph (nor strictly Deza, having diameter
$1$). So no Deza graph with a $2$-clique cover has $p\ne q$. If $p=q=m$, on the other hand, then
$p\,a=|E(H)|=q\,b$ forces $a=b=:t$, so $H$ is $t$-regular.
\end{proof}

\begin{remark}\label{rem:asym}
Theorem~\ref{thm:asym} settles the asymmetric case completely. Note that its proof uses only the
\emph{regularity} of $\GK$ together with the existence of the two-clique cover; the Deza property is
used only to conclude that an asymmetric two-clique cover is impossible, and no independence-number
argument is needed. The exhaustive search reported in
Table~\ref{tab:asym} is consistent with the theorem: for each of the pairs $(p,q)$ listed there, the
search found no Deza graph of diameter $2$; the regular configurations with $\dd_1=q$ and $\dd_2=p$
are exactly the complete ones, which are rejected at the diameter-$2$ stage and are not Deza graphs in
any case. So the table illustrates
Theorem~\ref{thm:asym} rather than proving it.
\end{remark}

\begin{corollary}[Classification for two-clique covers]\label{cor:2cliquecomplete}
Let $\GK$ be a Deza graph with a $2$-clique cover of parts of sizes $p$ and $q$. Then $p=q=:m$, and
exactly one of the
following holds.
\begin{enumerate}[label=\textup{(\alph*)}]
  \item $H$ is empty, that is, $t=0$. Then $\GK=2K_m$, which is regular of valency $m-1$
        with value set $\{0,m-2\}$ (for $m=1$ there are no pairs inside a clique, so the value set is the singleton $\{0\}$); here $m\ge3$, since for $m\le2$ only the single value $0$ occurs, so that $\GK$ is not a Deza graph. It is disconnected, so it is not strictly Deza.
  \item $H$ is $t$-regular with $1\le t\le m-1$. Put $c=m-t$. Then
        \[
          \ThL(\GK)=\bigl\{(m-2)+s:s\in\Sigma\cup\Sigma'\bigr\}\cup\{2t-2\},
          \qquad
          \Thm(\GK)=\{2t\},
        \]
        where the two values $2t-2<2t$ are distinct and both are attained; consequently
        \begin{gather*}
          \GK\ \text{is a Deza graph}
          \iff
          \Sigma\cup\Sigma'\subseteq\{t-c,\,t-c+2\},\\[2pt]
          \text{equivalently}\quad
          \bigl|\bigl((m-2)+(\Sigma\cup\Sigma')\bigr)\cup\{2t-2,2t\}\bigr|=2 .
        \end{gather*}
        Among the members of this family that \emph{are} Deza graphs, the two alternatives are
        separated by which of the two admissible values occurs in $\Sigma\cup\Sigma'$:
        \[
          \GK\ \text{is strongly regular}
          \iff
          \Sigma\cup\Sigma'=\{t-c\},
          \qquad
          \GK\ \text{is strictly Deza}
          \iff
          t-c+2\in\Sigma\cup\Sigma' .
        \]
\end{enumerate}
In case \textup{(b)}, $\Thm(\GK)=\{2t\}$ is always a singleton, so $\GK$ fails to be strongly regular
exactly when $\ThL(\GK)$ is not a singleton. For those members of case \textup{(b)} that \emph{are}
Deza graphs --- for which $\Sigma\cup\Sigma'\subseteq\{t-c,t-c+2\}$, so that
$\ThL(\GK)\subseteq\{2t-2,2t\}$ with $2t-2$ always attained --- this is in turn equivalent to
$2t=(m-2)+(t-c+2)$ occurring in $\ThL(\GK)$. Without the Deza hypothesis that last equivalence can
fail: for $m=7$, $t=4$ and $H$ the complement of the Fano incidence graph (its rows are the
complements of the seven Fano lines) we have $c=3$ and $\Sigma=\Sigma'=\{2\}$, which meets neither
$t-c=1$ nor $t-c+2=3$; correspondingly $\ThL(\GK)=\{6,7\}$, $\Thm(\GK)=\{8\}$, so $\GK$ is not
strongly regular and yet $2t=8\notin\ThL(\GK)$ --- indeed this $\GK$ is not even a Deza graph. Note
also that $|\Th(\GK)|=2$ characterises the \emph{Deza} property and not strictness:
the strongly regular members of the family, such as $(m,t)=(4,3)$, also have $|\Th(\GK)|=2$.
\end{corollary}

\begin{proof}
By Theorem~\ref{thm:asym} a $2$-clique cover of a Deza graph cannot have $p\ne q$: the only regular
configurations with $p\ne q$ are the complete graphs $K_{p+q}$, whose single common-neighbour value
$p+q-2$ shows that they are not Deza graphs. Hence $p=q=m$. Then $\GK$ is regular, so
\eqref{eq:asymdeg} forces $\dd_1\equiv a$ and $\dd_2\equiv b$ to be constant, and
\eqref{eq:asymdc} with $p=q=m$ gives $a=b=:t$; that is, $H$ is $t$-regular for some $0\le t\le m$. If
$t=0$ then $H$ is empty and $\GK=2K_m$, which is case \textup{(a)}: its values are $0$ and $m-2$, so
it is a Deza graph exactly for $m\ge3$, and it is disconnected. If $t=m$ then $H=K_{m,m}$ and
$\GK=K_{2m}=K_{p+q}$, which is not a Deza graph and so does not occur under the hypothesis. If
$1\le t\le m-1$ we are in case \textup{(b)}, where
Theorem~\ref{thm:valuesets} applies and gives the two displayed value sets
(Remark~\ref{rem:tboundary} explains why the restriction on $t$ is needed), hence also the Deza
criterion. For the final clause, $(m-2)+(t-c)=2t-2$ and $(m-2)+(t-c+2)=2t$, so the elements of
$(m-2)+(\Sigma\cup\Sigma')$ are precisely $2t-2$ and $2t$ according to which of $t-c$ and $t-c+2$ lies
in $\Sigma\cup\Sigma'$; the value $2t-2$ is in any case also contributed by the cross pairs, so it is
$\ThL$ that is or is not a singleton.
\end{proof}

\begin{remark}[The union $\Sigma\cup\Sigma'$ is indispensable here]\label{rem:sigmaonly}
It is tempting to replace $\Sigma\cup\Sigma'$ in Corollary~\ref{cor:2cliquecomplete} by $\Sigma$ alone,
on the grounds that $\Sigma$ already records the intersections of the rows of $H$ and the rows
determine $H$. This is not a harmless simplification: the two versions give different answers. Let
$m=8$, $t=4$ and let $H$ have rows
\begin{gather*}
  00001111,\quad 00110011,\quad 00111100,\quad 01010101,\\
  10101010,\quad 11000011,\quad 11001100,\quad 11110000
\end{gather*}
as subsets of $[8]$, read from left to right. Then $H$ is $4$-regular on $8+8$ vertices and
$\Sigma=\{0,2\}=\{2t-m,2t-m+2\}$, so the $\Sigma$-only criterion would declare $\GK(8,4,H)$ strictly
Deza. But $\Sigma'=\{1,2,3\}$, and $\GK(8,4,H)$ is not even a Deza graph: it is regular of valency
$11$ with $\ThL=\{6,7,8,9\}$ and $\Thm=\{8\}$. Hence the $\Sigma$-only form of the criterion is
\emph{false}, not merely unproved. Note that this $H$ has the same $\Sigma$ as the strictly
Deza witness $(8,4)$ of Table~\ref{tab:witnesses}, and the same valency $11$; what separates the two
graphs is $\Sigma'$.
\end{remark}

\section{Open problems}\label{sec:open}

\begin{proposition}[A strictly Deza graph with clique cover number $3$]\label{prop:wagner}
Let $F$ be the Wagner graph on $\mathbb{Z}/8\mathbb{Z}$, in which $i$ is adjacent to $i\pm1$ and to
$i+4$, and let $\GK=\overline F$. Then $\GK$ is strictly Deza with parameters $(8,4,2,1)$, lies in
$\cG$, and has clique cover number exactly $3$, the cover being $\{0,2,5\}$, $\{1,3,6\}$, $\{4,7\}$.
\end{proposition}

\begin{proof}
$F$ is $3$-regular, so $\GK$ is $4$-regular on $8$ vertices. Checking the $28$ pairs gives
$\ThL(\GK)=\{1,2\}$ and $\Thm(\GK)=\{2\}$, so $\GK$ is strictly Deza with parameters $(8,4,2,1)$.
The colouring $0,1,0,1,2,0,1,2$ is proper on $F$, and $F$ is triangle-free, so $F$ is $3$-colourable
and triangle-free; by \eqref{eq:G-description} this gives $\GK\in\cG$. Finally, cliques of $\GK$ are
independent sets of $F$, so a $2$-clique cover of $\GK$ would be a proper $2$-colouring of $F$; but $F$
contains the $5$-cycle $0$--$1$--$2$--$3$--$4$--$0$ and is therefore not bipartite. Hence the clique
cover number is at least $3$, and the displayed partition consists of independent sets of $F$, so it is
exactly $3$.
\end{proof}

\begin{remark}\label{rem:wagner-known}
The graph $\GK$ of Proposition~\ref{prop:wagner} is not new: $F$ is the M\"obius ladder and
$\GK$ is isomorphic to the Cayley graph $\mathrm{Cay}(\mathbb Z_8;\{\pm1,\pm2\})$ --- the isomorphism is
multiplication by $3$ in $\mathbb Z_8$, which carries the connection set $\{\pm2,\pm3\}$ of $\GK$ to
$\{\pm1,\pm2\}$ --- and it is one of the four non-$2$-extendable
Deza graphs of even order and diameter $2$ listed in \cite[Fig.~1]{MiklavicSparl2015}; the introduction
of \cite{KutnarMarusicMiklavicSparl2019} names it as the only such graph besides $C_4$, the Petersen
graph and $K_{2,2,2}$. It is also one of
the three strictly Deza graphs on $8$ vertices of the enumeration \cite{GoryainovPanasenkoShalaginov2021},
the other two being the rows $(m,t)=(4,1)$ and $(4,2)$ of Table~\ref{tab:main}. What is used here is its
membership in $\cG$: $F$ is $3$-colourable and triangle-free, so $\GK$ is a prime graph of a finite
solvable group. For extendability results on Deza graphs of diameter $2$ see
\cite{MiklavicSparl2015,KutnarMarusicMiklavicSparl2019}.
\end{remark}

Here we leave two open problems for further study.

\begin{problem}[Three-clique covers]\label{prob:beyond}
By \eqref{eq:G-description} every graph in $\cG$ has a $3$-clique cover, while by
Theorem~\ref{thm:2clique-in-G} a \emph{$2$}-clique cover forces membership in $\cG$.
Proposition~\ref{prop:wagner} shows that strictly Deza graphs in $\cG$ with clique cover number $3$
exist; all the infinite families of this paper, and those rook graphs that lie in $\cG$
(Proposition~\ref{prop:uniform}(iv)), have a $2$-clique cover. Determine all strictly Deza graphs in $\cG$ whose \emph{minimum} clique cover
has three parts, and decide in particular whether infinitely many of them exist.
\end{problem}

\begin{problem}[Groups of Lie type]\label{prob:lie}
Determine whether a finite simple group of Lie type can have a strictly Deza prime graph. We are not
aware of a proof in either direction. A connection with the present framework may exist: $(0,2)$-graphs
arise in the study of the cohomology of the unipotent radical of groups of Lie type with simply laced
diagram \cite{BrouwerChastkofsky2013}. A Zsigmondy-type argument on element orders appears to be the
natural approach, but we have not been able to make it precise. Here ``strictly Deza'' is used in the
sense of Definition~\ref{def:deza}: a prime graph of diameter $1$ (that is, a complete graph) or of
infinite diameter is excluded by definition.
\end{problem}


\clearpage
\appendix
\section{Numerical tables}\label{app:tables}

\begin{table}[htbp]
\centering
\small
\caption{Witness bipartite graphs $H$ for the strictly Deza cases of Table~\ref{tab:main} and of
Problem~\ref{prob:allm}.}
\label{tab:witnesses}
\begin{tabular}{@{}l >{\raggedright\arraybackslash}p{11.8cm}@{}}
\toprule
$(m,t)$ & $H$ \\
\midrule
$(4,1)$ & $x_i\leftrightarrow y_i$ for $i=0,\dots,3$ \\
$(4,2)$ & $x_0,x_1\leftrightarrow y_0,y_1$;\quad $x_2,x_3\leftrightarrow y_2,y_3$ \\
$(6,4)$ & $x_0,x_1\leftrightarrow y_0,y_1,y_2,y_3$;\quad
          $x_2,x_3\leftrightarrow y_0,y_1,y_4,y_5$;\quad
          $x_4,x_5\leftrightarrow y_2,y_3,y_4,y_5$ \\
$(7,3)$ & here the displayed sets \emph{are} the rows $R_i$ of $H$, i.e.\ the seven blocks of a\
          $2$-$(7,3,1)$ design (a Fano plane):\
          $x_0\leftrightarrow y_0,y_1,y_2$;\ $x_1\leftrightarrow y_0,y_3,y_4$;\
          $x_2\leftrightarrow y_0,y_5,y_6$;\ $x_3\leftrightarrow y_1,y_3,y_5$;\
          $x_4\leftrightarrow y_1,y_4,y_6$;\ $x_5\leftrightarrow y_2,y_3,y_6$;\
          $x_6\leftrightarrow y_2,y_4,y_5$ \\
$(8,4)$ & $00001111$, $00110011$, $01010101$, $01101001$,\ $10010110$, $10101010$,\
          $11001100$, $11110000$ \\
$(8,5)$ & $00011111$, $00101111$, $01001111$, $10001111$,\ $11110001$, $11110010$,\
          $11110100$, $11111000$ \\
$(11,6)$ & the eleven rows $00000111111$, $00111000111$, $01011011001$, $01101101010$,\
          $01110110100$, $10011101100$, $10101110001$, $10110011010$, $11001010110$,\
          $11010100011$, $11100001101$ \\
$(12,7)$ & the twelve rows $000001111111$, $001110001111$, $010110110011$, $011011010101$,\
          $011101101001$, $100111010110$, $101011100011$, $101100111010$, $110011111000$,\
          $110100011101$, $111001001110$, $111110100100$ \\
$(16,10)$ & the rows are the complements $R_i=\{0,\dots,15\}\setminus P_i$ of the blocks $P_i$ of a\
          $2$-$(16,6,2)$ biplane on the point set $\{0,\dots,15\}$ (note that this row uses $0$-based\
          labels, unlike $[m]=\{1,\dots,m\}$ elsewhere); the $P_i$ are $\{0,1,2,3,4,5\}$, $\{0,1,6,7,8,9\}$,\
          $\{2,3,6,7,10,11\}$, $\{4,5,8,9,10,11\}$,\
          $\{2,4,6,8,12,13\}$, $\{3,5,7,9,12,13\}$, $\{0,5,6,10,12,14\}$, $\{1,4,7,11,12,14\}$,\
          $\{1,3,8,10,13,14\}$, $\{0,2,9,11,13,14\}$, $\{1,2,9,10,12,15\}$, $\{0,3,8,11,12,15\}$,\
          $\{0,4,7,10,13,15\}$, $\{1,5,6,11,13,15\}$, $\{2,5,7,8,14,15\}$, $\{3,4,6,9,14,15\}$ \\
\bottomrule
\end{tabular}
\end{table}

In Table~\ref{tab:witnesses} the entries are to be read as follows. In the first four rows the cross
edges of $H$ are listed as adjacencies $x_i\leftrightarrow y_j$, so that equivalently the row $R_i$ is
the set of $y_j$ displayed; in the next four rows the rows of $H$ are given as binary masks over the
columns $y_0,\dots,y_{m-1}$, read from left to right; and in the last row the blocks $P_i$ of a
$2$-$(16,6,2)$ biplane are listed, the rows of $H$ being their complements
$R_i=\{0,\dots,15\}\setminus P_i$ (this row uses $0$-based labels). All nine rows were re-verified
independently: each $H$ is $t$-regular, the
resulting $\GK(m,t,H)$ is regular of valency $m-1+t$ (namely $4,5,9,9,11,12,16,18$ for the first
eight rows and $25$ for $(16,10)$, in agreement with Table~\ref{tab:main} in the four cases where the
row occurs there), and $\Sigma\cup\Sigma'$ satisfies \eqref{eq:criterion}.

\clearpage

\begin{table}[htbp]
\centering
\caption{Search-tree statistics for the exhaustive computation of Theorem~\ref{thm:smallm}.}
\label{tab:searchtree}
\begin{tabular}{@{}rrrrr@{}}
\toprule
$m$ & $t$ & search-tree nodes & graphs & Deza \\
\midrule
$2$ & $0$ & $3$ & $1$ & $0$ \\
$2$ & $1$ & $5$ & $2$ & $2$ \\
$2$ & $2$ & $3$ & $1$ & $0$ \\
$3$ & $0$ & $4$ & $1$ & $1$ \\
$3$ & $1$ & $16$ & $6$ & $0$ \\
$3$ & $2$ & $19$ & $6$ & $6$ \\
$3$ & $3$ & $4$ & $1$ & $0$ \\
$4$ & $0$ & $5$ & $1$ & $1$ \\
$4$ & $1$ & $65$ & $24$ & $24$ \\
$4$ & $2$ & $247$ & $90$ & $18$ \\
$4$ & $3$ & $109$ & $24$ & $24$ \\
$4$ & $4$ & $5$ & $1$ & $0$ \\
$5$ & $0$ & $6$ & $1$ & $1$ \\
$5$ & $1$ & $326$ & $120$ & $0$ \\
$5$ & $2$ & $5\,331$ & $2\,040$ & $0$ \\
$5$ & $3$ & $7\,471$ & $2\,040$ & $0$ \\
$5$ & $4$ & $901$ & $120$ & $120$ \\
$5$ & $5$ & $6$ & $1$ & $0$ \\
$6$ & $0$ & $7$ & $1$ & $1$ \\
$6$ & $1$ & $1\,957$ & $720$ & $0$ \\
$6$ & $2$ & $171\,631$ & $67\,950$ & $0$ \\
$6$ & $3$ & $967\,281$ & $297\,200$ & $0$ \\
$6$ & $4$ & $393\,361$ & $67\,950$ & $1\,350$ \\
$6$ & $5$ & $10\,051$ & $720$ & $720$ \\
$6$ & $6$ & $7$ & $1$ & $0$ \\
$7$ & $0$ & $8$ & $1$ & $1$ \\
$7$ & $1$ & $13\,700$ & $5\,040$ & $0$ \\
$7$ & $2$ & $7\,679\,743$ & $3\,110\,940$ & $0$ \\
$7$ & $3$ & $207\,870\,216$ & $68\,938\,800$ & $151\,200$ \\
$7$ & $4$ & $332\,053\,681$ & $68\,938\,800$ & $0$ \\
$7$ & $5$ & $32\,823\,526$ & $3\,110\,940$ & $0$ \\
$7$ & $6$ & $142\,297$ & $5\,040$ & $5\,040$ \\
$7$ & $7$ & $8$ & $1$ & $0$ \\
\midrule
\multicolumn{2}{@{}r}{total, $1\le t\le m-1$} & $582\,141\,934$ & $144\,548\,572$ & $158\,504$ \\
\multicolumn{2}{@{}r}{total, $0\le t\le m$} & $582\,142\,000$ & $144\,548\,584$ & $158\,509$ \\
\bottomrule
\end{tabular}
\end{table}

The columns of Table~\ref{tab:searchtree} are to be read as follows. The enumeration runs over row
masks with column-capacity pruning, a \emph{node} being a partial assignment of rows and a \emph{graph}
a completed assignment, i.e.\ one $t$-regular bipartite graph $H$ on $m+m$ labelled vertices. The last
column counts the labelled graphs $\GK(m,t,H)$ that are Deza graphs by definition, verified directly
from the adjacency matrix. The criterion \eqref{eq:criterion} presupposes $1\le t\le m-1$ and so does
not apply to the rows $t=0$ and $t=m$; those rows are listed only for completeness and are evaluated
directly. There $H$ is unique: for $t=0$ it is empty and $\GK=2K_m$, with both values attained exactly
when $m\ge3$ (hence the entry $1$ in the rows $(m,0)$ with $m\ge3$, and $0$ in $(m,t)=(2,0)$); for
$t=m$, $H=K_{m,m}$ and $\GK=K_{2m}$ has a single value, so the rows $(m,m)$ carry $0$. The ``graphs''
column was reproduced by an independent dynamic-programming count, and the whole table by a second
enumeration that uses no formula of this paper.

\clearpage

\begin{table}[htbp]
\centering
\caption{Rook graphs $\Lambda(r,p)$: value sets, from Proposition~\ref{prop:uniform}.}
\label{tab:uniform}
\begin{tabular}{@{}rrrrlll@{}}
\toprule
$r$ & $p$ & $v$ & $k$ & $\ThL$ & $\Thm$ & Deza? \\
\midrule
$2$ & $2$  & $4$  & $2$ & $\{0\}$   & $\{2\}$ & yes (SRG $C_4$) \\
$2$ & $3$  & $6$  & $3$ & $\{0,1\}$ & $\{2\}$ & no \\
$2$ & $4$  & $8$  & $4$ & $\{0,2\}$ & $\{2\}$ & \textbf{yes (strictly Deza, $=\GK_8$)} \\
$2$ & $5$  & $10$ & $5$ & $\{0,3\}$ & $\{2\}$ & no \\
$3$ & $3$  & $9$  & $4$ & $\{1\}$   & $\{2\}$ & yes (SRG, rook) \\
$3$ & $4$  & $12$ & $5$ & $\{1,2\}$ & $\{2\}$ & \textbf{yes (strictly Deza)} \\
$3$ & $5$  & $15$ & $6$ & $\{1,3\}$ & $\{2\}$ & no \\
$4$ & $4$  & $16$ & $6$ & $\{2\}$   & $\{2\}$ & no (SRG, $\lambda=\mu=2$: one value) \\
$4$ & $5$  & $20$ & $7$ & $\{2,3\}$ & $\{2\}$ & \textbf{yes (strictly Deza)} \\
$4$ & $6$  & $24$ & $8$ & $\{2,4\}$ & $\{2\}$ & \textbf{yes (strictly Deza)} \\
$5$ & $5$  & $25$ & $8$ & $\{3\}$   & $\{2\}$ & yes (SRG, rook) \\
$5$ & $6$  & $30$ & $9$ & $\{3,4\}$ & $\{2\}$ & no ($\Th=\{2,3,4\}$) \\
\bottomrule
\end{tabular}
\end{table}

\begin{table}[htbp]
\centering
\caption{Asymmetric two-clique covers: exhaustive search over all cross graphs $H$, under the
additional requirement that $\GK(p,q,H)$ be regular and of diameter $2$.}
\label{tab:asym}
\begin{tabular}{@{}rrrl@{}}
\toprule
$p$ & $q$ & number of $H$ & Deza graphs found \\
\midrule
$2$ & $3$ & $64$          & $0$ \\
$2$ & $4$ & $256$         & $0$ \\
$2$ & $5$ & $1\,024$      & $0$ \\
$2$ & $6$ & $4\,096$      & $0$ \\
$2$ & $7$ & $16\,384$     & $0$ \\
$3$ & $4$ & $4\,096$      & $0$ \\
$3$ & $5$ & $32\,768$     & $0$ \\
$4$ & $5$ & $1\,048\,576$ & $0$ \\
\bottomrule
\end{tabular}
\end{table}

The search reported in Table~\ref{tab:asym} enumerates \emph{all} bipartite graphs $H$ with parts of
the stated sizes, i.e.\ the number in the third column is $2^{pq}$, and then tests each resulting
$\GK(p,q,H)$ for regularity and diameter $2$ before testing the Deza condition. The complete graph
$\GK(p,q,K_{p,q})=K_{p+q}$ is rejected at the diameter-$2$ stage; by Theorem~\ref{thm:asym} the
complete configurations are the only regular ones with $p\ne q$, so no Deza graph arises there and the
table illustrates that theorem rather than proving it.


\end{document}